\documentclass{amsart}

\usepackage{amssymb,amscd,amsxtra,calc}
\usepackage{cmmib57}

\usepackage{amssymb}
\renewcommand{\baselinestretch}{1}

\newtheorem{theorem}{Theorem}[section]
\newtheorem{lemma}[theorem]{Lemma}
\newtheorem{proposition}[theorem]{Proposition}

\newtheorem{setup}[theorem]{}

\newtheorem{question}[theorem]{Question}

\newtheorem{_definition}[theorem]{Definition}

\newtheorem{_remark}[theorem]{Remark}
\newenvironment{remark}{\begin{_remark}\rm}{\end{_remark}}

\newtheorem{_difficulty}[theorem]{The difficulty}
\newenvironment{difficulty}{\begin{_difficulty}\rm}{\end{_difficulty}}

\newtheorem{_claim}[theorem]{Claim}
\newenvironment{claim}{\begin{_claim}\rm}{\end{_claim}}
\numberwithin{equation}{section}
\numberwithin{table}{section}
\numberwithin{figure}{section}

\newcommand{\C}{\mathord{\mathbb C}}

\newcommand{\D}{\Delta}

\newcommand{\N}{\mathord{\mathbb N}}
\renewcommand{\P}{\mathord{\mathbb P}}
\newcommand{\Q}{\mathord{\mathbb Q}}
\newcommand{\R}{\mathord{\mathbb R}}
\newcommand{\Z}{\mathord{\mathbb Z}}

\newcommand{\OO}{\mathord{\mathcal O}}

\newcommand{\NN}{\mathsf{N}}
\newcommand{\NEb}{\operatorname{\overline{\mathrm{NE}}}}
\newcommand{\Cont}{\operatorname{cont}}
\newcommand{\isom}{\simeq}

\newcommand{\Br}{{\text{\rm Br}}}

\newcommand{\CH}{{\text{\rm CH}}}

\newcommand{\diag}{{\text{\rm diag}}}

\newcommand{\Hom}{\text{\rm Hom}}
\newcommand{\id}{\text{\rm id}}

\newcommand{\NE}{\text{\rm NE}}
\newcommand{\Nef}{\text{\rm Nef}}

\newcommand{\Proj}{\text{\rm Proj}}
\newcommand{\rank}{\text{\rm rank}}
\newcommand{\Sing}{\text{\rm Sing}}

\newcommand{\Supp}{\text{\rm Supp}}

\newcommand{\Weil}{{\text{\rm Weil}}}

\newcommand{\ratmap}
{{\,\cdot\negmedspace\cdot\negmedspace\cdot\negmedspace\to\,}}

\begin{document}
\begin{large}
\title
[Polarized endomorphisms of uniruled varieties] {Polarized endomorphisms of uniruled varieties \\
(with Appendix by Y. Fujimoto and N. Nakayama) \\
\vskip 2pc
De-Qi Zhang}

%
%

%
%
%
%

\begin{abstract}
We show that polarized endomorphisms of rationally
connected threefolds with at worst terminal singularities are
equivariantly built up from those on $\Q$-Fano threefolds,
Gorenstein log del Pezzo surfaces and $\P^1$.
Similar results are obtained
for polarized endomorphisms of uniruled threefolds and fourfolds.
As a consequence, we show 
that every smooth Fano threefold
with a polarized endomorphism of degree $> 1$, is rational.
\end{abstract}
\subjclass{14E20, 14J45, 14E08, 32H50}
\keywords{polarized endomorphism, uniruled variety, rationality of variety}
\maketitle
\section{Introduction}
We work over the field $\C$ of complex numbers.
We study {\it polarized} endomorphisms
$f : X \to X$ of varieties $X$, i.e., those $f$ with
$f^*H \sim qH$ for some $q > 0$ and some
ample line bundle $H$.
Every surjective endomorphism of a projective variety of Picard number one, is polarized.
If $f = [F_0:F_1:\cdots:F_n] : \P^n \to \P^n$ is a surjective morphism and $X \subset \P^n$ a $f$-stable subvariety,
then $f^*H \sim qH$ and hence $f|X : X \to X$ is polarized; here $H \subset X$ is a hyperplane
and $q = \deg(F_i)$. If $A$ is an abelian variety and $m_A : A \to A$ the multiplication
map by an integer $m \ne 0$, then $m_A^*H \sim m^2 H$ and hence $m_A$ is polarized; here
$H = L + (-1)^*L$ with $L$ an ample divisor, or $H$ is any ample divisor with $(-1)^*H \sim H$.
One can also construct polarized endomorphisms on quotients of $\P^n$ or $A$.
So there are many examples of polarized endomorphisms $f$. See \cite{Zs} for the many conjectures on such $f$.

From the arithmetical point of view, given a polarized endomorphism $f : X \to X$
of degree $q^{\dim X}$ and
defined over $\Q$, one can define a unique height function $h_f : X(\overline{\Q}) \to \R$
such that $h_f(f(x)) = qf(x)$. Further, $x$ is $f$-preperiodic if and only if
$h_f(x) = 0$; see \cite[\S 4]{Zs} for more details.

In \cite{nz2}, it is proved that a normal variety $X$ with a
non-isomorphic polarized endomorphism $f$ either has only canonical singularities
with $K_X \sim_{\Q} 0$ (and further is a quotient of an abelian variety
when $\dim X \le 3$), or is uniruled so that $f$ descends to a polarized
endomorphism $f_Y$ of the non-uniruled base variety $Y$ (so $K_Y \sim_{\Q} 0$) of a specially chosen
maximal rationally connected fibration $X \ratmap Y$. By the induction on dimension
and since $Y$ has a dense set of $f_Y$-periodic points $y_0, y_1, \dots$ (cf. \cite[Theorem 5.1]{Fa}),
the study of polarized endomorphisms is then reduced to that of rationally connected varieties
$\Gamma_{y_i}$ as fibres of the graph $\Gamma = \Gamma(X/Y)$
(cf. \cite[Remark 4.3]{nz2}).

The study of non-isomorphic endomorphisms of singular varieties (like $\Gamma_{y_i}$ above)
is very important from the dynamics point of view, but is very hard even in dimension two
and especially for rational surfaces;
see \cite{Fv}, and \cite{ENS} (about 150 pages).

In this paper, we consider polarized endomorphisms of rationally connected varieties
(or more generally of uniruled varieties) of dimension $\ge 3$.
Theorem \ref{Th3-4} -- \ref{FinExt} below
and Theorems \ref{Th3-4a} -- \ref{FinExt1} in \S 3,
are our main results.

\begin{theorem}\label{Th3-4}
Let $X$ be a $\Q$-factorial $n$-fold, with $n \in \{3, 4\}$,
having only log terminal singularities and a polarized endomorphism $f$ of degree $q^n > 1$.
Let $X = X_0 \ratmap X_1 \cdots$ $\ratmap X_r$ be a composite of
divisorial contractions and flips.
Replacing $f$ by its positive power, we have:
\begin{itemize}
\item[(1)] The dominant rational maps
$g_i : X_i \ratmap X_i$ $(0 \le i \le r)$ $($with $g_0 = f)$
induced from $f$, are all holomorphic.
\item[(2)]
Let $\pi : X_r \to Y$ be an
extremal contraction
with $\dim Y \le 2$.
Then $g_r$ is polarized and it descends to a polarized endomorphism $h : Y \to Y$
of degree $q^{\dim Y}$ with
$\, \pi \circ g_r = h \circ \pi \, $.
\end{itemize}
\end{theorem}

The result above reduces the study of
$(X, f)$ to $(X_r, g_r)$ where the latter
is easier to be dealt with
since $X_r$ has a fibration structure
preserved by $g_r$.
{\it The existence of such a fibration $\pi: X_r \to Y$
is guaranteed} when $X$ is uniruled
by the recent development in MMP.
The relation between the two pairs is very close
becuase $f^{-1}$, as seen in Theorem \ref{Th3-4a}, 
preserves the maximal subset of $X$ where the birational map
$X \ratmap X_r$ is not holomorphic.

\begin{theorem}\label{Thev3}
Let $X$ be a $\Q$-factorial threefold having
only terminal singularities and a polarized endomorphism
of degree $q^3 > 1$.
Suppose that $X$ is rationally connected.
Then we have :
\begin{itemize}
\item[(1)]
There is an $s > 0$
such that ${{(f^s)}^*}_{| N^1(X)} = q^{s} \, \id$.
We then call such $f^s$ cohomologically a scalar.
\item[(2)]
Either $X$ is rational, or $-K_X$ is big.
\item[(3)]
There are only finitely many irreducible divisors $M_i \subset X$
with the Iitaka $D$-dimension $\kappa(X, M_i) = 0$.
\end{itemize}
\end{theorem}

Theorem \ref{Thev3} (3) above apparently does not hold 
for $X = S \times \P^1$, where $S$ is a rational surface
with infinitely many $(-1)$-curves and hence $S$ has no endomorphisms
of degree $> 1$ by \cite[Proposition 10]{Ny02};
the blowup of nine general points of $\P^2$
is such $S$ as observed by Nagata.

Theorem \ref{Thev3} (1) above strengthens (in our situation) Serre's result \cite{Se}
on a conjecture of Weil (in the projective case): ({\bf Serre}) If $f$ is a polarized
endomorphism of degree $q^{\dim X} > 1$ of a smooth variety $X$ then every
eigenvalue of $f^* | N^1(X)$ has the same modulus $q$.

The proof of Theorem \ref{Corrat3} below is done
without using the classification of smooth Fano threefolds.
This result has been reproved in \cite{ratcon} where $f$ is
assumed to be only of degree $> 1$ but not necessarily polarized.
%

\begin{theorem}\label{Corrat3}
Let $X$ be a smooth Fano threefold with
a polarized endomorphism $f$ of degree $> 1$.
Then $X$ is rational.
\end{theorem}

A klt $\Q$-Fano variety has only finitely many extremal rays.
A similar phenomenon occurs in the quasi-polarized case (cf. \ref{conv}).

\begin{theorem} \label{FinExt}
Let $X$ be a $\Q$-factorial rationally connected threefold
having only Gorenstein terminal singularities and a quasi-polarized endomorphism of degree $> 1$.
Then $X$ has only finitely many $K_X$-negative extremal rays.
\end{theorem}


%

The claim in the abstract about the building blocks
of polarized endomorphisms, is justified by the remark below.

\begin{remark}

$ $

\par
(1) The $Y$ in Theorem \ref{Th3-4} is $\Q$-factorial and has at worst log terminal singularities;
see \cite{ZDA}.

(2) Suppose that the $X$ in Theorem \ref{Th3-4} is rationally connected.
Then $Y$ is also rationally connected.
Suppose further that $X$ has at worst terminal singularities and $(\dim X, \ \dim Y) = (3, 2)$.
Then $Y$ has at worst Du Val singularities by \cite[Theorem 1.2.7]{MP}.
So there is a composition $Y \to \hat{Y}$ of divisorial contractions
and an extremal contraction ${\hat Y} \to B$
such that either $\dim B = 0$ and ${\hat Y}$ is a Du Val del Pezzo
surface of Picard number $1$, or $\dim B = 1$
and ${\hat Y} \to B \cong \P^1$ is a $\P^1$-fibration with all fibres irreducible.
After replacing $f$ by its power, $h$ descends to polarized endomorphisms
$\hat{h} : {\hat Y} \to {\hat Y}$, and
$k : B \to B$ (of degree $q^{\dim B}$); see Theorems \ref{Thev2a}.

(3) By \cite[Theorem 5.1]{Fa}, there are dense subsets $Y_0 \subset Y$
(for the $Y$ in Theorem \ref{Th3-4})
and $B_0 \subset B$ (when $\dim B = 1$)
such that for every $y \in Y_0$ (resp. $b \in B_0$)
and for some $r(y) > 0$ (resp. $r(b) > 0$),
$g^{r(y)} | W_y$ (resp. $\hat{h}^{r(b)} | {\hat Y}_b$)
is a well-defined polarized endomorphism of the Fano fibre.
\end{remark}

\begin{difficulty}
In Theorem \ref{Th3-4}, if $X \to X_1$ is a divisorial contraction, one can descend
a polarized endomorphism $f$ on $X$ to an one on $X_1$, but the latter may not be polarized any more
because the pushfoward of a nef divisor may not be nef in dimension $\ge 3$
(the first difficulty).
If $X \ratmap X_1$ is a flip, then in order to descend $f$ on $X$ to some holomorphic $f_1$ on $X_1$,
one has to show that a power of $f$ preserves the centre of the flipping contraction
(the second difficulty).
The second difficulty is taken care by Lemma \ref{Mfin} where the polarizedness is
essentially used. 

As pointed out by the referee, a key argument in the proof of
Theorem \ref{Th3-4} (2) is to show that a power of $f$ is cohomologically a scalar
unless $Y$ is a surface with torsion $K_Y$ (this case will not happen when $X$
is rationally connected); see Lemma \ref{pi}.
\end{difficulty}

The question below is the generalization of Theorem \ref{Corrat3} and
the famous conjecture: every smooth Fano $n$-fold of {\it Picard number one}
with a non-isomorphic surjective endomorophism, is $\P^n$
(for its affirmative
solution when $n = 3$, see  Amerik-Rovinsky-Van de Ven \cite{ARV} and Hwang-Mok \cite{HM}).

\begin{question}\label{question}
Let $X$ be a smooth Fano $n$-fold with a non-isomorphic polarized endomorphism.
Is $X$ rational ?
\end{question}

\begin{remark}
A recent preprint of Koll\'ar and Xu \cite{KX} showed that one can descend the endomorphism
$\P^n \to \P^n$ ($[X_0, \dots, X_n] \to [X_0^m, \dots, X_n^m]$; $m \ge 2$)
to some quotient $X := \P^n/G$ (with $G$ finite) so that $X$ has only terminal singularities
but $X$ is irrational, invoking a famous prime power order 
group action of David Saltman on Noether's problem.
Thus one cannot remove the smoothness assumption in Theorem \ref{Corrat3} and Question \ref{question}.

However, we will show in Theorem \ref{Thrat3} that
every rationally connected $\Q$-factorial projective threefold $X$
with only terminal singularities, is rational, provided that
$X$ has a non-isomorphic polarized endomorphism
and an extremal contraction $X \to Y$ with $\dim Y \in \{1, 2\}$.
The terminal singularity assumption there is used to deduce the Gorenstein-ness
of $Y$ (when $\dim Y = 2$),
making use of \cite[Theorem 1.2.7]{MP}. 
 
As pointed out by the referee, it would be interesting if one
could determine
whether the `terminal singularity' assumption
can further be weakened to the `log canonical singularity'
in order to deduce the rationality as above.

See also \cite{ratcon} for the generalization
of Theorem \ref{Thrat3} to non-polarized endomorphisms.
\end{remark}

For the recent development
on endomorphisms of algebraic varieties,
we refer to Amerik-Rovinsky-Van de Ven \cite{ARV}, Fujimoto-Nakayama \cite{FN},
Hwang-Mok \cite{HM},
Hwang-Nakayama \cite{HN}, S. -W. Zhang \cite{Zs}, as well as \cite{NZ}, \cite{ICCM}.
\par \vskip 1pc \noindent
{\bf Acknowledgement}

I would like to thank N. Nakayama for informing me
about his paper \cite{ENS} and critical comments on the draft of the paper. I thank
Y. Fujimoto and N. Nakayama
for patiently listening to my talks at RIMS and the constructive suggestions,
the algebraic geometry group at RIMS, Kyoto University : S. Mori,
S. Mukai, N. Nakayama, \dots for the stimulating atmosphere
and the warm hospitality
in the second half of the year 2007,
and the referee for the constructive comments.

I also like to thank A. Fujiki of Osaka Univ., S. Kondo of Nagoya Univ.,
Y. Kawamata of Univ. of Tokyo, and A. Moriwaki of Kyoto Univ. for
providing me with the opportunity of talks.

\par
This project is supported by an Academic Research Fund of
NUS.
\section{Preliminary results}
\begin{setup}\label{conv}
{\bf Conventions}
\end{setup}

{\it Every endomorphism in this paper is assumed to be surjective.}

For a projective variety $X$, an endomorphism $f : X \to X$
is {\it polarized} or {\it polarized by $H$} (resp. {\it quasi-polarized}
or {\it quasi-polarized by $H$})
if $f^*H \sim_{\Q} qH$ for some $q > 0$ and some
ample (resp. nef and big) line bundle $H$.
If $f$ is polarized or quasi-polarized then so
is its induced endomorphism on the normalization of $X$.

On a projective variety $X$, denote by
$N^1(X)$ (resp. $N_1(X)$) the usual
$\R$-vector space of $\R$-Cartier $\R$-divisors
(resp. $1$-cycles with coefficients in $\R$)
modulo numerical equivalence, in terms of the
perfect pairing
$N^1(X) \times N_1(X) \to \R$.
The Picard number $\rho(X)$ equals $\dim_{\R} N^1(X) = \dim_{\R} N_1(X)$.
The {\it nef cone} $\Nef(X)$ is the closure in $N^1(X)$ of the ample cone,
and is dual to the closed cone $\overline{\NE}(X) \subset N_1(X)$
generated by effective $1$-cycles (Kleiman's ampleness criterion).

Denote by $S(X)$ the set of $\Q$-Cartier prime divisors $G$
with $G_{|G}$ non-pseudo-effective; see \cite[II, \S 5]{ZDA}
for the relevant material.

For a normal projective surface $S$, a Weil divisor
is {\it numerically equivalent to zero} if so is its Mumford pullback
to a smooth model of $S$. Denote by $\Weil(S)$ the set of $\R$-divisors (divisor = Weil divisor)
modulo this numerical equivalence. We can also define the intersection of
two Weil divisors by Mumford-pulling back them to a smooth model and then taking
the usual intersection.

A Weil divisor is {\it nef} if its intersection with
every curve is non-negative.
A Weil divisor $D$ on a normal projective variety is {\it big}
if $D \sim_{\Q} A + E$ for an ample line bundle $A$ and an effective
Weil $\R$-divisor $E$ (see \cite[II, 3.15, 3.16]{ZDA}).

Let $f : X \to X$ be an endomorphism and $\sigma_V : V \to X$ and $\sigma_Y : X \to Y$
morphisms. We say that $f$ {\it lifts} to an endomorphism $f_V : V \to V$ if
$f \circ \sigma_V = \sigma_V \circ f_V$;
$f$ {\it descends} to an endomorphism $f_Y$ if
$\sigma_Y \circ f = f_Y \circ \sigma_Y$.

A normal projective variety $X$ is $Q$-{\it abelian} in the sense of \cite{nz2}
if $X = A/G$ with $A$ an abelian variety and $G$ a finite group acting
freely in codimension $1$, or equivalently $X$ has an abelian variety
as an \'etale in codimension 1 cover.

For a normal projective variety $X$, we refer to \cite{KMM}
or \cite{KM} for the definition of $\Q$-{\it factoriality} and {\it terminal singularity}
or {\it log terminal singularity}.
An {\it extremal contraction} $X \to Y$ is always assumed to be $K_X$-negative.

We do not distinguish a Cartier divisor with its corresponding line bundle.

\begin{lemma}\label{ev}
Let $X$ be a normal projective $n$-fold and $f : X \to X$
an endomorphism such that $f^*H \equiv qH$ for some $q > 0$ and a nef and big
line bundle $H$. Then we have:

\begin{enumerate}
\item[(1)]
There is a nef and big line bundle $H'$ such that
$H' \equiv H$ and $f^*H' \sim_{\Q} qH'$. So $f$ is quasi-polarized.
Further, $\deg(f) = q^n$.
\item[(2)]
Every eigenvalue of $f^* | N^1(X)$ has modulus $q$.
\item[(3)]
Suppose that $\sigma: X \to Y$ is a fibred space (with connected fibres) and $f$ descends to
an endomorphism $h : Y \to Y$. Then $\deg(h) = q^{\dim Y}$.
Every eigenvalue of $h^* | N^1(Y)$ has modulus $q$.
\end{enumerate}
\end{lemma}

\begin{proof}
(1) and (2) are just \cite[Lemmas 2.1 and 2.3]{nz2}.

Set $d := \deg(h)$ and $\dim Y = k$.
Then $f^*X_y \equiv dX_y$ for a general fibre $X_y$ over
$y \in Y$. Now (3) follows
from the fact that $\sigma^*N^1(Y)$ is a $f^*$-stable subspace of $N^1(X)$
and the calculation:
$$q^n H^{n-k} . X_y = f^*H^{n-k} . f^*X_y = q^{n-k} d H^{n-k} . X_y > 0.$$
\end{proof}

\begin{setup} {\bf Pullback of cycles}
\end{setup}
We will consider pullbacks of cycles by finite surjective morphisms.
Let $X$ be a normal projective variety.
We define a numerical equivalence $\equiv$ for cycles
in the Chow group $\CH_r(X)$ of $r$-cycles modulo rational equivalence.
An $r$-cycle is called {\it numerically equivalent to zero}, denoted as $C \equiv 0$,
if $H_1 \dots H_r . C = 0$ for all Cartier divisors $H_i$.

If $C$ is a nonzero effective $r$-cycle then $C$ is not numerically equivalent to
zero since $H^r . C > 0$ for an ample line bundle $H$.
Denote by $[C]$ the equivalence class of all $r$-cycles numerically
equivalent to $C$.
Denote by $N_r(X)$ the set $\{[C] \, ; \, C$ is an $r$-cycle
with coefficients in $\R \}$.
The usual product of an $r$-cycle
with $s$ line bundles naturally extends to
$$N^1(X) \times \cdots \times N^1(X) \times N_r(X) \,\, \longrightarrow \,\,
N_{r-s}(X).$$

Let $f : X \to X$ be a surjective endomorphism of degree $d$,
so $f$ is a finite morphism.
For an $r$-dimensional subvariety $C$, write $f^{-1} C = \cup_i C_i$
and define $f^*[C] := \sum_i e_i [C_i]$ with $e_i > 0$ chosen such that
$\sum_i e_i \delta_i = d$ for $\delta_i := \deg(C_i/C)$.
Then
$$f_*f^*[C] = d[C].$$
If $C, C_i$ are not in $\Sing X$, then for the usual
$f^*$-pullback $f^*C$ of the cycle $C$, we have $[f^*C] = f^*[C]$ by having the right choice of $e_i$.
By the linearity of the intersection form,
we can linearly extend the definition to $f^*[C]$ for an arbitrary $r$-cycle $C$.
Then the usual projection formula gives
$$f^*L_1 \dots  f^*L_r . f^*[C] = \deg(f) (L_1 \dots L_r . C).$$
Note that $f^* : N^1(X) \to N^1(X)$ is an isomorphism.
With this, $[C] \to f^*[C]$ (or simply $f^*C$ by the abuse of notation) 
gives a well defined map
$$f^*: N_r(X) \,\, \longrightarrow N_r(X).$$
The projection formula above implies the following in $N_{r-s}(X)$
$$f^*(L_1 \dots L_s . C) \equiv f^*L_1 \dots f^*L_s . f^*C.$$

\begin{lemma}\label{deg}
Let $X$ be a normal projective $n$-fold and $f : X \to X$ an endomorphism
of degree $q^n$ for some $q > 0$. Suppose that every eigenvalue
of $f^* | N^1(X)$ has modulus $q$.  Then we have:
\begin{itemize}
\item[(1)]
If $D$ is an $r$-cycle such that $0 \ne [D] \in N_r(X)$
and $f^*D \equiv a D$. Then $|a| = q^{n-r}$.
%

\item[(2)]
If $S$ is a $k$-dimensional subvariety of $X$ with $f^{-1}(S) = S$ as set, then
$f^*S \equiv q^{n-k} S$ and $\deg(f : S \to S) = q^{k}$.

\item[(3)]
Suppose the $S$ in $(2)$ is a surface.
Then there is a Cartier $\R$-divisor $M$ on $X$
such that $M_S := M_{|S}$ is a nonzero element in $\Nef(S)$ and $f_{\,\,\,|S}^*M_S \equiv qM_S$
in $N^1(S)$.

\item[(4)]
If $\rho(X) \le 2$, then $(f^2)^*|N^1(X) = q^2 \ \id$.
\end{itemize}
\end{lemma}

\begin{proof}
(4) We may assume that $(f^2)^*E_i \equiv a_i E_i$ for the extremal
rays $E_i$ ($1 \le i \le \rho(X)$) in $\Nef(X)$.
Thus $a_i = |a_i| = q^2$ by the assumption, done!

(2) follows from (1) and our definition of pullback.

(1) Choose a basis $L_1, \dots, L_{\rho}$ with $\rho = \rho(X)$
such that $f^*|N^1(X)$ is lower triangular.
So $f^*L_i = q u(i) L_i +$ lower term with $|u(i)| = 1$.
Since $[D] \ne 0$, for some $s > 0$,
the cycle $L_{s} . D$ is not numerically equivalent to zero.
We choose $s$ to be minimal.
Now
$$f^*(L_{s} . D) \equiv f^*L_s . f^*D =
(q u(s) L_s + \text{lower term}) . aD = aqu(s) (L_s . D).$$
Similarly, we can show that
$C : = L_s . L_{s_1} \dots L_{s_{r-2}} . D \in N_1(X)$ is
not numerically equivalent to zero, and
$f^*C \equiv b C$ with
$$b = aq^{r-1} \prod_{i=0}^{r-2} u(s_i), \hskip 1pc (s_0 := s).$$
Since $N_1(X)$ is dual to $N^1(X)$, the eigenvalue $b$ of $f^*|N_1(X)$
satisfies $|b| = q^{n-1}$. So $|a| = q^{n-r}$ as claimed.

(3) Let $N^1(X)_{|S} \subseteq N^1(S)$ (resp. $\Nef(X)_{|S} \subseteq \Nef(S)$)
be the image of $\iota^* : N^1(X) \to N^1(S)$ (resp. of the restriction of this $\iota^*$
to $\Nef(X)$)
with $\iota: S \to X$ the closed embedding.
Let $\overline{N}$ be the closure of $\Nef(X)_{|S}$ in $N^1(S)$.
Then $\overline{N}$ spans the subspace $N^1(X)_{|S}$ of $N^1(S)$.
Let $\lambda$ be the spectral radius of $f^*|\overline{N}$.
By the generalized Perron-Frobinius theorem in \cite{Bi},
$f^*(M_S) \equiv \lambda (M_S)$ for a nonzero nef divisor
$M_S := M_{|S}$ in $\overline{N}$ (with $M$ a Cartier $\R$-divisor on $X$).
Write $M|S = a_t L_t | S +$ lower term, with $t$ the smallest (and $a_t \ne 0$).
Then
$$\lambda a_t L_t | S + \text{lower term} =
\lambda M|S = f^*(M|S) = a_t q u(t) L_t|S + \text{lower term}.$$
By the minimality of $t$, we have $\lambda a_t = a_t q u(t)$
and $\lambda = |\lambda| = q$.
\end{proof}

\begin{lemma}\label{ev'}
Let $X$ be a normal projective surface and $f : X \to X$ an endomorphism
of degree $q^2 > 1$. Suppose that $f^*M \equiv qM$ for a nonzero nef
Weil divisor.
Then every eigenvalue of $f^*|\Weil(X)$ has modulus $q$.
\end{lemma}

\begin{proof}
Let $\lambda$ be the spectral radius of $f^* | \Weil(X)$.
Then $f^*L \equiv \lambda L$ for a nonzero nef $\R$-divisor $L$.
Now $q^2 L . M = f^*L . f^*M = \lambda q L . M$. So either
$L . M > 0$ and $\lambda = q$, or $L . M = 0$. In the latter case,
$M \equiv c L$ by the Hodge index theorem (on a resolution of $X$)
and again we have $\lambda = q$.

Similarly, let $\mu$ be the spectral radius of $(f^*)^{-1} | \Weil(X)$
so that $(f^*)^{-1} H \equiv \mu H$ for a nonzero nef $\R$-divisor $H$.
Then $f^*H \equiv \mu^{-1} H$. By the argument above, we have $\mu^{-1} = q$.
The lemma follows.
\end{proof}

Here is an easy polarizedness criterion for ruled normal surfaces.

\begin{lemma}\label{ruled}
Let $X$ be a normal projective surface and $X \to B$ a $\P^1$-fibration.
Suppose that $f : X \to X$ is an endomorphism of degree $q^2 > 1$
and $f^*H \equiv qH$ for a nonzero nef $\R$-divisor $H$.
Then there is an $s > 0$ such that $(f^s)^*|\Weil(X) = q^s \id$.
So $f$ is polarized.
\end{lemma}

\begin{proof}
Note that a basis of $\Weil(X)$ consists of some negative curves $C_1, \dots, C_r$
in fibres,
a general fibre and a multiple section. Contract
$C_i$' to get a Moishezon normal surface $Y$ with
$\Weil(Y) = \R E_1 + \R E_2$ for two extremal rays
$\R_{\ge 0} E_i$ of the cone $\overline{\NE}(X)$.
By \cite[Proposition 10]{Ny02} or as in the proof of Lemma \ref{S},
replacing $f$ by its power, we may assume that $f^{-1}(C_i) = C_i$ for all $i$.

So $f$ descends to an endomorphism $f_Y : Y \to Y$
and we may assume that $f^*E_i \equiv e_i E_i$ for some $e_i > 0$
after replacing $f$ by $f^2$.

Write $f^*C_i = a_iC_i$ with $a_i > 0$.
Then $f^* | \Weil(X) = \diag[a_1, \dots, a_r, e_1, e_2]$
with respect to the basis: $C_1, \dots, C_r$ and the pullbacks of $E_1, E_2$.
Now the first assertion follows from Lemma \ref{ev'}
while the second follows from the first as in Note 1 of Theorem \ref{Thev2a}.
This proves the lemma.
\end{proof}

Nakayama's \cite[Example 4.8]{ENS} (ver. Jan 2008) produces many
examples of polarized $f$ on abelian surfaces which are not scalar.
The result below shows that this happens only on abelian surfaces and their quotients.

\begin{theorem}\label{Thev2a}
Let $X$ be a normal projective surface.
Suppose that $f : X \to X$ is an endomorphism
such that $f^*P \equiv qP$ for some $q > 1$ and
some big Weil $\Q$-divisor $P$. Then we have:
\begin{itemize}
\item[(1)]
$f$ is polarized of degree $q^2$.
\item[(2)]
There is an $s > 0$ such that $(f^s)^*|\Weil(X) = q^s \, \id$
unless $X$ is $Q$-abelian with $\rank \Weil(X) \in \{3, 4\}$.
\end{itemize}
\end{theorem}

\begin{proof}
Let $P = P' + N'$ be the Zariski decomposition. Then $P'$ is a nef and big
Weil $\Q$-divisor. The uniqueness of such decomposition and $f^*P \equiv qP$
imply $f^*P' \equiv qP'$ and $f^*N' \equiv qN'$.
Replacing $P$ by $P'$, we may assume that $P$ is already a nef and big
Weil $\R$-divisor. So $\deg(f) = (f^*P)^2/P^2 = q^2$.

\par \vskip 1pc
{\bf Note 1.} If $(f^s)^* H' \equiv q^s H'$
for an ample line bundle $H'$ on $X$ then $f$ is polarized.
Indeed, If we set $H:= \sum_{i=0}^{s-1} (f^i)^*H/q^i$,
then $H$ is an ample $\Q$-divisor with $f^*H \equiv qH$, and we apply Lemma \ref{ev}.

\par \vskip 1pc
{\bf Claim 1}. $ $
\begin{itemize}
\item[(1)]
Every eigenvalue of $f^*|\Weil(X)$ has modulus $q$.
\item[(2)]
If $(f^s)^* | \Weil(X)$
is scalar for some $s > 0$, then
it is $q^s \, \id$.
\end{itemize}
\par \vskip 1pc

Claim 1(1) follows from Lemma \ref{ev'} while Claim 1(2) follows from (1).

\par \vskip 1pc
Claim 2 below is from Claim 1 and the proof of Lemma \ref{deg} (4).

\par \vskip 1pc
{\bf Claim 2}. If $\rho := \dim_{\R} \Weil(X) \le 2$, then
$(f^2)^* | \Weil(X) = q^2 \, \id$.

\par \vskip 1pc

\par \vskip 1pc
By \cite[Proposition 10]{Ny02} or as in the proof of Lemma \ref{S}, 
the set $S'(X)$ of negative curves on $X$
is finite and $f^{-1}$ induces a bijection of $S'(X)$.
We may assume that $f|S'(X) = \id$ after replacing $f$ by its power.
Let $X \to Y$ be the composition of contractions of negative curves $C_1, \dots, C_r$
(with $r$ maximum)
intersecting the canonical divisor negatively.
Then $Y$ is a relatively minimal Moishezon normal surface in the sense of \cite{Sa}.
$f$ descends to an endomorphism $f_Y : Y \to Y$.

\par \vskip 1pc
{\bf Case(1)} $K_Y$ is not pseudo-effective.
Then either $\rank \Weil(Y) = 2$ and there is a $\P^1$-fibration $Y \to B$, or
$\Weil(Y) = \R[-K_Y]$ with $-K_Y$ numerically ample;
see \cite[Theorem 3.2]{Sa}. With $f$ replaced by its square,
we may assume that $f_Y^* | \Weil(Y) = q \, \id$ (use Claim 1, and see the proof of
Lemma \ref{deg} (4)).
Thus $f^* | \Weil(X) = q \, \id$
with respect to the basis
consisting of $C_1, \dots, C_r$ and the pullback of a basis of
$\Weil(Y)$; see Claim 1.
So the theorem is true in this case.

\par \vskip 1pc
{\bf Case(2)} $K_Y$ is pseudo-effective (and hence nef by the minimality).
So $K_X$ is also pseudo-effective.
It is well known then that the ramification divisor $R_f = 0$ and hence
$f$ is \'etale in codimension $1$. Further, $K_X = f^*K_X$ and hence
$K_X^2 = 0$ since $\deg(f) > 1$.
If $C \in S'(X)$ is a negative curve on $X$ then
$f^*C = qC$ by Claim 1, and because of the extra assumption $f | S'(X) = \id$,
$f$ is ramified along $C$. Thus $S'(X) = \emptyset$.
So $X = Y$ and $K_X$ is nef. Also $P$ is numerically ample.
The proof is completed by:

\par \vskip 1pc
{\bf Claim 3}. $X$ is $Q$-abelian. So $\rank \Weil(X) \le 4$,
$X$ is $\Q$-factorial,
and $f$ is polarized by $P$ which is $\Q$-Cartier.

\par \vskip 1pc
Since $q^2P . K_X = f^*P . f^*K_X = qP . K_X$, we have $P . K_X = 0$.
The Hodge index theorem (applied to a resolution of $X$)
implies that $K_X \equiv 0$ in $\Weil(X)$.
Thus the claim follows from \cite[Theorem 7.1.1]{ENS}.
\end{proof}

\begin{lemma}\label{degbir}
Let $X$ be a normal projective $n$-fold and $f : X \to X$
a quasi-polarized endomorphism of degree $q^n > 0$.
Then we have:

\begin{enumerate}
\item[(1)]
Suppose that $V \to X$ is a birational morphism and $f$ lifts
to an endomorphism $f_V : V \to V$. Then $f_V$ is also quasi-polarized.
\item[(2)]
Let $X \ratmap W$ be a birational map with $W$ being $\Q$-factorial,
such that the dominant
rational map $f_W : W \ratmap W$ induced from $f$, is holomorphic.
Then $f_W^*H_W \sim_{\Q} qH_W$ for some big line bundle $H_W$
and every eigenvalue of $f_W^*|N^1(W)$ has modulus $q$.
\end{enumerate}
\end{lemma}

\begin{proof}
By the definition, there is a line bundle $H$ on $X$
such that $f^*H \sim_{\Q} qH$.
(1) holds because $f_V$ is quasi-polarized by the pullback $H_V$ of $H$.

(2) Let $V$ be the normalization of the graph $\Gamma_{X/W}$. Then
$f$ lifts to a quasi-polarized endomorphism $f_V$ of $V$.
For the first assertion, we take $H_W$ to be (a multiple of) the direct image of $H_V$
(consider pullback to $V$ of $H_W$ and use Lemma \ref{ev} (2) and the argument in Note 1 of Theorem
\ref{Thev2a}).
Since $N^1(W)$ can be regarded as a subspace of $N^1(V)$
with the action $f_W^*$ and $f_V^*$ compactible,
the second follows from Lemma \ref{ev}.
\end{proof}

\begin{lemma}\label{S}
Let $V$ and $X$ be normal projective $n$-folds with $X$ being $\Q$-factorial, and
$\tau : V \ratmap X$ a birational map. Suppose
an endomorphism $f : X \to X$ of degree $> 1$, lifts to a quasi-polarized endomorphism
$f_V : V \to V$. Then the set $S(X)$ of prime divisors $D$ on $X$
with $D_{|D}$ not pseudo-effective, is a finite set. Further,
$f^{-1}(S(X)) = S(X)$, so $f^r|S(X) = \id$ for some $r > 0$.
\end{lemma}

\begin{proof}
Replacing $V$ by the normalization of the graph of $\tau : V \ratmap X$
and using Lemma \ref{degbir}, we may assume that
$\tau$ is already holomorphic.
By the assumption, there is a nef and big line bundle $H$
such that $f_V^*H \sim qH$ and hence $\deg(f) = \deg(f_V) = q^n > 1$.
Note that $f^*$ and $f_* = q^n (f^*)^{-1}$ are automorphisms on both $N^1(X)$
and $N_1(X)$.

{\bf Step 1}. If $D \in S(X)$ then $D':= f(D) \in S(X)$. Indeed, $f^*D' \equiv cD$
with $c > 0$ because $f_*(f^*D')$ is parallel to $f_*D$.
Since $f^*({D'}_{|D'}) \equiv cD_{|D}$ is not
pseudo-effective, $D' \in S(X)$.

{\bf Step 2}. If $D' := f(D) \in S(X)$ then $D \in S(X)$. This is because
$f^*D' \equiv cD$ as in Step 1 and hence
$cD_{|D} \equiv f^* ({D'}_{|D'})$ is not pseudo-effective.

{\bf Step 3}. If $f(D_1) = D' = f(D_2)$ for $D_1 \in S(X)$, then $D_1 = D_2$.
Indeed, $f_*D_1 \equiv e f_*D_2$ for some $e > 0$. So $D_1 \equiv eD_2$.
Since $eD_{2|D_1} \equiv D_{1|D_1}$ is not pseudo-effective, $D_1 = D_2$.

It follows then

{\bf Step 4}. $f^{-1}(S(X)) = S(X)$, and $f$ and $f^{-1}$ act bijectively on $S(X)$.

{\bf Step 5}. Let $(H^{n-1})^{\perp}$ be the set of prime divisors $F$
with $F . H^{n-1} = 0$. Then it is a finite set.
Indeed, writing $H = A + E$ with $A$ an ample Cartier $\Q$-divisor and $E$
an effective Cartier $\Q$-divisor, then the set above is contained in the support of $E$.

{\bf Step 6}. There is a finite set $\Sigma$, such that
$f^{c(D)}(D) \in \Sigma$ with some $c(D) \ge 0$ for every $D \in S(X)$.
This will imply the lemma (see \cite[Proposition 10]{Ny02}). 
We take $\Sigma$ to be the union of the set of prime divisors in
$\Sing X$ 
and the ramification divisor
$R_{f}$ of $f$, and the set
of prime divisors on $X$
whose strict transform on $V$ is in $(H^{n-1})^{\perp}$.

To finish Step 6, we only need to consider those $D \in S(X)$ where $D_i := f^{i-1}(D)$ is
not in $\Sigma$ for all $i \ge 1$.
Write $f^*D_{i+1} = a_iD_i$ with $a_i \in \Z_{>0}$.
Let $D_i' \subset V$ be the strict transform of $D_i$.
Then $f_V^*D_{i+1}' \equiv a_i D_i'$ in $N_{n-1}(V)$.
So
$$q^n H^{n-1} . D_{i+1}' = f_V^*H^{n-1} . f_V^*D_{i+1}' = q^{n-1}a_i H^{n-1} . D_i',$$
$$1 \le H^{n-1} . D_{i+1}' = \frac{a_i}{q} \cdots \frac{a_1}{q}H^{n-1} . D_1'.$$
Thus $a_{i_0} \ge q$ for infinitely many $i_0$. So $D_{i_0}$ is in $R_f$
and hence in $\Sigma$. This completes Step 6 and also the proof of the lemma.
\end{proof}

\begin{lemma}\label{Mfin}
Let $V$ and $X$ be projective $n$-folds, $\tau : V \to X$
a birational morphism, $\D = \D_X \subset X$ a Zariski-closed subset
and $f : X \to X$ an endomorphism of degree $q^n > 1$.
Assume the four conditions below:
\begin{itemize}
\item[(1)]
$f$ lifts to an endomorphism $f_V : V \to V$ quasi-polarized
by a nef and big line bundle $H$ so that $f^*H \sim qH$.
\item[(2)]
$f^{-1}(\D(i)) = \D(i)$ for every irreducible component $\D(i)$ of $\D$
(but we only need $f^{-1}(\D) = \D$ in the proof).
\item[(3)]
$\tau : V \to X$ is isomorphic over $X \setminus \D$.
\item[(4)]
For every subvariety $Z \subset V$ not contained in $\tau^{-1}(\D)$,
the restriction $H_{|Z}$ is nef and big
(and hence $\deg(f|Z : Z \to Z) = q^{\dim Z}$).
\end{itemize}
Let $A \subset X$ be a positive-dimensional subvariety such that $f^{-j}f^j(A) = A$
for all $j \ge 0$. Then either $M(A) := \{f^i(A) \, | \, i \ge 0\}$
is a finite set, or $f^{i_0}(A) \subseteq \D$ for some $i_0$
(and hence for all $i \ge i_0$).
\end{lemma}

\begin{proof}
We shall prove by induction on the codimension of $A$ in $X$.

Set $k := \dim A$, $A_1 := A$ and $A_i := f^{i-1}(A)$ ($i \ge 1$).
Denote by $\Sigma$ or $\Sigma(V, X, \D, f)$
the set of prime divisors in $\D$, $\Sing X$ and
the ramification divisor $R_f$ of $f$.
This $\Sigma$ is a finite set.

\par \vskip 1pc
{\bf Claim 1}. $A_i$ is contained in the union $U(\Sigma)$
of prime divisors in $\Sigma$ for infinitely many $i$;
so if $\dim A = \dim X - 1$, our $M(A)$ is finite and the lemma holds.
\par \vskip 1pc

Suppose the contrary that Claim 1 is false. Replacing $A$ by some $A_{i_0}$,
we may assume that $A_j$ is not contained in $U(\Sigma)$ for all $j \ge 1$.
Set $b_j := \deg(f: A_j \to A_{j+1})$.
Write $f^*A_{j+1} = a_j A_j$ as cycles with $a_j = q^n/b_j \in \Z_{> 0}$ now.
Let $A_j' \subset V$ be the strict transform of $A_j$.
Now $f_V^*A_{j+1}' = a_j A_j'$ as cycles, and
$$q^n H^k . A_{j+1}' = f_V^*H^k . f_V^*A_{j+1}' = q^k a_j H^k . A_j',$$
$$1 \le H^k . A_{j+1}' = \frac{a_j}{q^{n-k}} \cdots \frac{a_1}{q^{n-k}}H^k . A_1'.$$
Thus $a_{j_0} \ge q^{n-k}$ for infinitely many $j_0$.
So $A_{j_0}$ is contained in $R_f$ and hence also in $U(\Sigma)$ for infinitely many $j_0$.
Thus Claim 1 is true.

\par \vskip 1pc
We may assume that $|M(A)| = \infty$ and $k \le n - 2$.
Let $B$ be the Zariski-closure of the union of those $A_{i_0}$ contained in $U(\Sigma)$.
Then $\dim B \in \{k+1, \dots, n-1\}$, and $f^{-j}f^j(B) = B$ for all $j \ge 0$.
Choose $r \ge 1$ such that $B' := f^r(B), f(B'), f^2(B'), \dots$ all
have the same number of irreducible components.
Let $X_1$ be an irreducible component of $B'$ of maximal dimension.
Then $\dim X_1 \in \{k+1, \dots, n-1\}$ and $f^{-j}f^j(X_1) = X_1$ for all $j \ge 0$.
Note also that $X_1$ contains infinitely many $A_{i_1}$.
If $f^j(X_1) \subseteq \D$ for some $j \ge 0$, then $A_{i_1+j} \subseteq \D$
and we are done. Thus we may assume that
$\D \cap f^j(X_1) \subset f^j(X_1)$ for all $j \ge 0$ and hence $M(X_1) < \infty$
by the inductive assumption with codimension. We may assume that
$f^{-1}(X_1) = X_1$,
after replacing $f$ with its power and $X_1$ with its image of some $f^j$.

Let $V_1 \subset V$ be the strict transform of $X_1$. Then
all four conditions in the lemma are satisfied by $(V_1, H|V_1, X_1, \D|X_1, f|X_1, A_{i_1})$.
Since the codimension of $A_{i_1}$ in $X_1$ is smaller than that of $A$ in $X$,
by the induction, either $M(A_{i_1})$ and hence $M(A)$ are finite
or $A_{j_0} \subseteq \D | X_1 \subseteq \D$ for some $j_0$.
This completes the proof of the lemma.
\end{proof}

\begin{lemma}\label{ext}
Let $X$ be a projective variety and $f : X \to X$ a surjective
endomorphism. Let $R_C := \R_{\ge 0}[C] \subset \overline{\NE}(X)$
be an extremal ray (not necessarily $K_X$-negative).
Then we have:
\begin{itemize}
\item[(1)]
$R_{f(C)}$ is an extremal ray.
\item[(2)]
If $f(C_1) = C$, then $R_{C_1}$ is an extremal ray.
\item[(3)]
Denote by $\Sigma_C$ the set of curves whose classes are in $R_C$.
Then $f(\Sigma_C) = \Sigma_{f(C)}$.
\item[(4)]
If $R_{C_1}$ is extremal then $\Sigma_{C_1} = f^{-1}(\Sigma_{f(C_1)})
:= \{D \, | \, f(D) \in \Sigma_{f(C_1)}\}$.
\end{itemize}
\end{lemma}

\begin{proof}
Note that $f^* : N^1(X) \to N^1(X)$ and $f_* : N_1(X) \to N_1(X)$
are isomorphisms.

(1) Suppose $z_1 + z_2 \equiv f_*C$ for $z_i \in \overline{\NE}(X)$.
Write $z_i = f_*z_i'$ for $z_i' \in \overline{\NE}(X)$. Then
$f_*(z_1'+z_2' -C) \equiv 0$ and hence $z_1' + z_2' \equiv C$.
Thus $z_i' \equiv a_iC$ for some $a_i \ge 0$ by the assumption on $C$, whence $z_i
= f_*z_i' \equiv a_i f_*C \in R_{f(C)}$.

(2) $\sim$ (4) are also easy.
\end{proof}

\begin{lemma}\label{-Kext}
Let $X$ be a normal projective variety with at worst log terminal singularities,
and $f : X \to X$ an endomorphism. Suppose that $R_{C_i} = \R_{\ge 0}[C_i]$ ($i = 1, 2$),
with $C_2 = f(C_1)$,
are $K_X$-negative extremal rays and $\pi_i : X \to Y_i$ the corresponding contractions.
Then there is a finite surjective morphism $h : Y_1 \to Y_2$ such that $\pi_2 \circ f = h \circ \pi_1$.
\end{lemma}

\begin{proof}
Let $X \to Y \overset{h}\to Y_2$ be the Stein factorization of
$\pi_2 \circ f : X \to X \to Y_2$. By Lemma \ref{ext}, the map $X \to Y$
is just $\pi_1 : X \to Y_1$.
\end{proof}

The result below is crucial and used in proving Theorem \ref{Th3-4a}.
It was first proved by the author when $\dim Y \le 2$ or $\rho(Y) \le 2$,
and has been extended and simplified by Fujimoto and Nakayama to the current form below.
See Appendix for its proof.

\begin{theorem}\label{descend}
Let \( X \) be a normal projective variety defined over an algebraically
closed field of characteristic zero such that \( X \)
has only log-terminal singularities.
Let \( R \subset \NEb(X) \) be an extremal ray such that \( K_{X}R < 0 \)
and the associated contraction morphism \( \Cont_{R} \) is a fibration to
a lower-dimensional variety.
Then, for any surjective endomorphism \( f \colon X \to X \),
there exists a positive integer \( k \) such that
\( (f^{k})_{*}(R) = R \) for the automorphism
\( (f^{k})_{*} \colon \NN_{1}(X) \xrightarrow{\isom} \NN_{1}(X)\)
induced from the iteration \( f^{k} = f \circ \cdots \circ f \).
\end{theorem}

\section{Proof of Theorems}

In this section we prove the theorems in the Introduction and three theorems below.
Theorem \ref{Th3-4a} below includes Theorem \ref{Th3-4} as
a special case, while Theorem \ref{FinExt1} implies \ref{FinExt} 
because a result of Benveniste says that a Gorenstein terminal threefold has no flips.
We note:

\begin{remark} \label{RemS3}
All $X_i$, $Y$ in Theorem \ref{Th3-4a} are again $\Q$-factorial
and have at worst log terminal singularities by MMP (see e.g. \cite{ZDA}).
\end{remark}

\begin{theorem}\label{Th3-4a}
Let $X$ be a $\Q$-factorial $n$-fold, with $n \in \{3, 4\}$,
having only log terminal singularities and a polarized endomorphism $f$ of degree $q^n > 1$.
Let $X = X_0 \ratmap X_1 \cdots$ $\ratmap X_r$ be a composite of
$K$-negative divisorial contractions and flips.
Replacing $f$ by its positive power, {\rm (I)} and {\rm (II)} hold:
\begin{itemize}
\item[(I)] The dominant rational maps
$g_i : X_i \ratmap X_i$ $(0 \le i \le r)$ $($with $g_0 = f)$ induced from $f$,
are all holomorphic. Further, $g_i^{-1}$ preserves each irreducible
component of the exceptional locus of $X_i \to X_{i+1}$ $($when it is
divisorial$)$ or of the flipping contraction $X_i \to Z_i$ $($when
$X_i \ratmap X_{i+1} = X_i^{+}$ is a flip$)$.
\item[(II)]
Let $\pi : W = X_r \to Y$ be the contraction of a $K_W$-negative extremal ray
$\R_{\ge 0}[C]$,
with $\dim Y \le n-1$.
Then $g := g_r$ descends to a surjective endomorphism $h : Y \to Y$
of degree $q^{\dim Y}$ such that
$$\pi \circ g = h \circ \pi.$$
For all $0 \le i \le r$,
all eigenvalues of $g_i^*|N^1(X_i)$ and $h^*| N^1(Y)$ are of modulus $q$;
there are big line bundles $H_{X_i}$ and $H_Y$ satisfying
$$
g_i^*H_{X_i} \sim q H_{X_i}, \hskip 1pc h^*H_Y \sim qH_Y.$$
Suppose further that either $\dim Y \le 2$ or $\rho(Y) = 1$. Then
$H_{W}$ and $H_Y$ can be chosen to be ample and $g$ and $h$ are
polarized.
\end{itemize}
\end{theorem}

The contraction $\pi$ below exists by the MMP for threefolds.

\begin{theorem}\label{Thrat3}
Let $X$ be a $\Q$-factorial rationally connected threefold
having at worst terminal singularities and a polarized endomorphism of degree $> 1$.
Let $X \ratmap W$ be a composite of $K$-negative divisorial contractions and flips, and
$\pi: W \to Y$ an extremal contraction of non-birational type.
Suppose either $\dim Y \ge 1$, or $\dim Y = 0$ and $W$ is smooth.
Then $X$ is rational.
\end{theorem}

\begin{theorem} \label{FinExt1}
Let $X$ be a $\Q$-factorial rationally connected threefold
having only terminal singularities.  
Suppose either $X$ has a quasi-polarized endomorphism of degree $> 1$,
or the set $S(X)$ as in $\ref{conv}$ is finite.
Then $X$ has only finitely many $K_X$-negative extremal rays which are not of flip type.
\end{theorem}

We start with some preparations for the proof of Theorem \ref{Th3-4a}.

\begin{proposition}\label{pol2}
Let $X$ be a $\Q$-factorial $n$-fold with $n \in \{3, 4\}$, having at worst
log terminal singularities and a polarized endomorphism
$f: X \to X$ of degree $q^n > 1$.
Let $X = X_0 \ratmap X_1 \cdots \ratmap X_r$ be a composite of
$K$-negative divisorial contractions and flips.
Suppose that for each $0 \le j \le r$,
the dominant rational map $f_j : X_j \ratmap X_j$
induced from $f$, is holomorphic and $f_j^{-1}$ preserves
each irreducible component of the exceptional locus of
$X_j \to X_{j+1}$ (when it is divisorial) or of the flipping contraction $X_j \to Y_j$
(when $X_j \ratmap X_{j+1} = X_j^{+}$ is a flip).
Let $S'$ be a surface on some $X_i$ with $(f_i^v)(S') = S'$
for some $v > 0$.
Then the endomorphism $f_S: S \to S$ induced from $f_i^v|S'$,
is polarized of degree $q^{2v}$. Here
$S$ is the normalization of $S'$.
\end{proposition}

\begin{proof}
We may assume that $v = 1$ after replacing $f$ by its power;
see Note 1 of Theorem 2.7.
By the assumption, $f^*H_{X} \sim qH_X$ for a very ample line bundle $H_X$,
and $\deg(f) = q^n$.
By Lemmas \ref{degbir} and \ref{deg}, $\deg(f_S : S \to S) = q^2$.
To show the polarizedness of $f_S$,
we only need to show the assertion of the existence of
a big Weil divisor as an eigenvector
of $f_S^*$; see Theorem \ref{Thev2a}.

We shall prove this assertion by ascending induction on the index $i$ of $X_i$.
When $X_i = X$, $S$ is polarized by the pullback of $H_X$
via the morphism $S \to S' \subset X$.

If $X_{i-1} \to X_{i}$ is birational over $S'$ with
$S_{i-1}' \subset X_{i-1}$ the strict transform of $S'$
and $S_{i-1}$ the normalization of $S_{i-1}'$, then
the polarizedness of $S_{i-1}$ (by the inductive assumption)
gives rise to a big Weil divisor $P_S$ on $S$
with $f_S^* P_S \equiv q P_S$ (using Lemma \ref{ev'} and the proof of
Lemma \ref{degbir}). We are done.

Thus, we have only to consider the two cases below
(where $n = 4$).

{\bf Case(1)} $X_{i-1} \to X_i$ is a divisorial contraction
so that $S'$ is the image of a prime divisor $Z'$ on $X_{i-1}$
(being necessarily the support of the whole exceptional divisor $X_{i-1} \to X_i$).
By the assumption, $f_{i-1}^{-1}(Z') = Z'$ and hence $f^{-1}(Z_X') = Z_X'$
where $Z_X' \subset X$ is the (birational) strict transform of $Z'$.
The normalization $Z$ of $Z_X'$ has an endomorphism $f_Z$ (induced from $f|Z_X'$)
polarized by $H_Z$ (the pullback of $H_X$) so that $f_Z^*H_Z \sim qH_Z$. $Z' \to S'$ induces
$\sigma: Z \to S$ (with general fibre $\P^1$) so that $f_S$ is the descent of $f_Z$.
By \cite[the proof of Proposition 4.17]{IntS},
the intersection sheaf $H_S := I_{Z/S}(H_Z, H_Z)$ is an integral Weil
divisor satisfying $f_S^* H_S \sim q H_S$.
Further, $H_S = (\sigma|H_Z)_*(H_{Z|H_Z})$ and hence is big
by the ampleness of $H_Z$. We are done again.

{\bf Case(2)} $X_{i-1} \ratmap X_i = X_{i-1}^+$ is a flip and $S'$ is an irreducible
component of the exceptional locus of the flipping contraction $X_i \to Y_{i-1}$.
We have $f_i^{-1}(S') = S'$ by the assumption on the flipping contraction
$X_{i-1} \to Y_{i-1}$.
Note that the assumption of Lemma \ref{deg} is satisfied by
$(X_i, f_i)$ (see Lemma \ref{degbir}). In particular,
$f_i^* M|S' \equiv q M|S'$ for a nonzero nef Cartier $\R$-divisor
$M |S'$ in $N^1(X_i)|S' \subset N^1(S')$. We divide into two subcases.

Case(2a) $S'$ is mapped to a curve $B'$ on $Y_{i-1}$.
Then we have an induced map $S \to B$ with general fibre $\P^1$.
Here $B$ the normalization of $B'$.
Thus $f_S$ is polarized by Lemmas \ref{ruled} and \ref{deg}.

Case(2b) $S'$ is mapped to a point on $Y_{i-1}$.
Note that $\rho(X_i/Y_{i-1}) = 1$ since $\rho(X_{i-1}/Y_{i-1}) = 1$
and $\rho(X_{i-1}) = \rho(X_i)$.
So for any ample Cartier divisor $A$ on $X_i$, there is a $b \ne 0$
such that $A - bM$ is the pullback of some divisor by $X_i \to Y_{i-1}$.
Thus $A|S' \equiv bM|S'$ in $N^1(S')$. Hence
$f_i^* A|S' \equiv qA|S'$ in $N^1(S')$.
Thus $f_S$ is polarized by an ample line bundle
$A_S$ (the pullback of $A|S'$).
\end{proof}

I thank N. Nakayama for suggesting the proof below.

\begin{lemma}\label{flip1}
Let $X$ be a $\Q$-factorial projective variety with at worst
log terminal singularities, $f : X \to X$ a surjective endomorphism,
and $X \ratmap X^+$ a flip with $\pi: X \to Y$ the corresponding
flipping contraction of an extremal ray $R_C := \R_{\ge 0}[C]$.
Suppose that $R_{f(C)} = R_C$. Then the dominant rational map
$f^+ : X^+ \ratmap X^+$ induced from $f$, is holomorphic.
Both $f$ and $f^+$ descend to one and the same endomorphism
of $Y$.
\end{lemma}

\begin{proof}
We note that
$$X = \Proj \oplus_{m \ge 0} \OO_Y(-mK_Y), \hskip 1pc
X^+ = \Proj \oplus_{m \ge 0} \OO_Y(mK_Y)$$
and there is a natural birational morphism
$\pi^+ : X^+ \to Y$. By the assumption and Lemma \ref{-Kext},
$f : X_1 = X  \to X_2 = X$ descends to an endomorphism $h: Y_1 = Y \to Y_2 = Y$
with $\pi_2 \circ f = h \circ \pi_1$. Here $\pi_i : X_i \to Y_i$
are identical to $\pi : X \to Y$.
Set $Z := X_2^+ \times_{Y_2} Y_1$. Then the projection $Z \to Y_1$
is a small birational morphism with $\rho(Z/Y_1) = 1$, and
it is identical to either $X_1 \to Y_1$
or $X_1^+ = X^+ \to Y_1$, noting that $-K_X$ and $K_{X^+}$ are relatively
ample over $Y$. Now we have only to consider and rule out the case $Z = X_1$.
Set $W := X_2^+ \times_{Y_2} X_2$. Since the composite $X_1 = Z \to X_2^+ \to Y_2$
is identical to that of $Z \to Y_1 \to Y_2$ and hence to that of
$X_1 \to X_2 \to Y_2$, there is a morphism $\sigma: X_1 \to W$
such that $X_1 = Z \to X_2^+$ factors as $X_1 \to W \to X_2^+$,
and $X_1 \to X_2$ factors as $X_1 \to W \to X_2$.
So the projection $W \to X_2$ is birational (because so is $X_2^+ \to Y_2$)
and finite (because so is $X_1 \to X_2$), whence it is an isomorphism.
Thus the birational map $X_2 \to X_2^+$ is a well defined morphism
as the composition of $X_2 \to W \to X_2^+$. This is absurd.
Therefore, $Z = X_1^+$ and the lemma is true.
\end{proof}

\begin{lemma}\label{div}
With the hypotheses and notation in Lemma $\ref{Mfin}$, assume further
that $X$ is $\Q$-factorial with at worst log terminal singularities and
$\sigma: X \to X_1$ is a divisorial contraction of an extremal ray $\R_{\ge 0}[\ell]$
with $E$ the exceptional locus (necessarily an irreducible divisor).
Then we have:
\begin{itemize}
\item[(1)]
There is an $s > 0$ such that
$(f^s)^{-1}(E) = E$.
\item[(2)]
The dominant rational map $g : X_1 \ratmap X_1$ induced from $f^s$,
is holomorphic, after $s$ is replaced by a larger one.
\item[(3)]
Let $\D_1 \subset X_1$ be the image of $\D \cup E$.
Then $g^{-1}(\D_1) = \D_1$.
\item[(4)]
Let $V_1$ be the normalization of the
graph of $V \ratmap X_1$, and $H_1 \subset V_1$ the pullback of $H$
on $V$. Then $g$ lifts to
an endomorphism $g_1 : V_1 \to V_1$ such that
$(V_1 \supset H_1, \, g_1, \, X_1 \supset \D_1, g)$
satisfies all four conditions in Lemma $\ref{Mfin}$.
\end{itemize}
\end{lemma}

\begin{proof}
(1) follows from Lemma \ref{S}
since $E \in S(X)$, while (3) and (4) follow from (2).
Now (2) follows from the proof of Theorem \ref{descend}
applied to $N^1(X) | E \subset N^1(E)$
and the extremal curve $\ell$ in the closed cone
of curves on $E$ (dual to the cone
$\Nef(X)|E$).
\end{proof}

\begin{lemma}\label{flip2}
With the hypotheses and notation in Lemma $\ref{Mfin}$, assume further:
\begin{itemize}
\item[(1)]
If $T' \subset X$ is a surface with $f^t(T') = T'$ for some $t > 0$,
then the endomorphism of the normalization $T$ of $T'$ induced from ${f^t}_{|T'}$,
is polarized.
\item[(2)]
$\dim \D \le 2$.
\item[(3)]
$X$ has at worst log terminal singularities and $X \ratmap X^+$ is a flip with $\pi : X \to Y$ the corresponding
flipping contraction of an extremal ray $R_C := \R_{\ge 0}[C]$.
\item[(4)]
The union $U_C$ of curves in the set $\Sigma_C$ in Lemma $\ref{ext}$ is of dimension $\le 2$.
\end{itemize}
Then we have:
\begin{itemize}
\item[(1)]
There is an $s > 0$ such that $R_{f^s(C)} = R_C$ and
$(f^s)^{-1}(U_C(i)) = U_C(i)$ for every irreducible component
$U_C(i)$ of $U_C$.
\item[(2)]
The dominant rational map $g : X^+ \ratmap X^+$ induced from $f^s$,
is holomorphic.
\item[(3)]
Let $\D^+ = \D(X^+) \subset X^+$ be the set consisting of the exceptional locus of the flipping
contraction $\pi^+ : X^+ \to Y$ (i.e., $(\pi^+)^{-1}(\pi(U_C))$)
and the total transform of $\D \subset X$. Then $g^{-1}(\D^+(i)) = \D^+(i)$
for every irreducible component $\D^+(i)$ of $\D^+$.
\item[(4)]
Let $V^+$ be the normalization of the
graph of $V \ratmap X^+$, and $H^+ \subset V^+$ the pullback of $H$
on $V$. Then $g$ lifts to
an endomorphism $g_{V^+} : V^+ \to V^+$ such that
$(V^+ \supset H^+, \, g_{V^+}, \, X^+ \supset \D^+, g)$
satisfies all four conditions in Lemma $\ref{Mfin}$.
\end{itemize}
\end{lemma}

\begin{proof}
Note that the assertion(2) follows from (1) and Lemma \ref{flip1}, while
(3) and (4) follow from (1) and (2).
It remains to prove (1). By Lemma \ref{ext},
we have only to show that $f^{u}(C)$ and $f^v(C)$ (and hence
$f^{u-v}(C)$ and $C$) are parallel for some $u > v$.

By Lemma \ref{ext}, $f^{-j}f^j(U_C) = U_C$ for all $j \ge 0$.
Choose $r' \ge 0$ such that $U':= f^{r'}(U_C), f(U'), f^2(U'), \dots$
all have the same number of irreducible components. Then $f^{-j}f^j(U'(k)) = U'(k)$
for every irreducible component $U'(k)$ of $U'$.
By Lemma \ref{Mfin}, either $M(U'(k))$ is finite
and $S' := f^{j_1}(U'(k)) = f^{j_2}(U'(k))$
for some $j_2 > j_1 > 1$, or $f^{j_1}(U'(k))$ is contained in
an irreducible component $\D(1)$ of $\D$
for infinitely many $j_1$.
We divide into two cases.

\vskip 1pc
Case(1) $\dim U'(k) = 2$. Since $\dim \D(1) \le 2$
we may assume that $M(U'(k))$ is always finite
and $(f^m)^{-1}(S') = S'$ for $m = j_2 - j_1$.
Take a $2$-dimensional irreducible component $S$ of $U_C$
such that $f^{r}(S) = S'$, where $r := r'+j_1$.
Note that $f^{-m}$ permutes irreducible components of
$f^{-r}(S')$. So some $f^{-t}$ with $t \in m \N$, stabilizes
all of these components. Especially, $f^{\pm t}(S) = S$.
Replacing $f$ by $f^t$,
we may assume that $f^{\pm}(S) = S$. We may also assume that $C \subset S$.
If the flipping contraction $\pi : X \to Y$
maps $S$ to a point $P$, then $f(C)$ is parallel to $C$ because $\pi(f(C)) = P$,
so (1) is true.
Suppose $\pi$ induces a fibration $S \to B$ onto a curve. Let
$\widetilde{S} \to S$ be the normalization. Then $f$ induces
a finite morphism $\widetilde{f}: \widetilde{S} \to \widetilde{S}$
which is polarized by our assumption, so $\widetilde{f}^* | \Weil(\widetilde{S})
= q \, \id$ after replacing $f$ by its power (see Lemmas \ref{ruled}, \ref{deg} and \ref{degbir}).
Thus $f(C)$ is parallel to $C$. Hence (1) is true in Case(1).

\par \vskip 1pc
Case(2) $\dim U'(k) = 1$. We only need to consider the situation where
$f^{j_1}(U'(k)) \subset \D(1)$ and $\dim \D(1) = 2$.
Relabel $f^{r'+j_1}(C)$ as $C$, we have $C \subset S := \Delta(1)$.
By the hypotheses, $f^{\pm}(S) = S$.
Set $C_v:= f^v(C)$.
By the choice of $r'$, we have $f^{-j}f^j(C) = C$ for all $j \ge 0$.
Let $\widetilde{S} \to S$ be the normalization and $\Theta \subset \widetilde{S}$
the union of the conductor and the ramification divisor $R_h$ of
the finite morphism $h : \widetilde{S} \to \widetilde{S}$ induced from $f$.
If $C_{v}$ has preimage in $\Theta$ for infinitely many $v$ then $C_{v}$ and $C_{v'}$
(and hence $C_{v - v'}$ and $C$) are parallel for some $v > v'$ because
$\Theta$ has only finitely many components, so (1) is true.
Thus we may assume that no $C_{v}$ is contained in $\Theta$ for all $v \ge 0$.
Let $D_v \subset \widetilde{S}$
be the birational preimage of $C_{v}$. Then $h^{-j}h^j(D_v) = D_v$ for all $j \ge 0$.
The extra assmuption implies $h^*D_{v+1} = D_v$. By Lemmas \ref{deg} and \ref{degbir},
we have $\deg(h) = q^{2}$.
Now $q^{2} D_{v+1} . D_{w+1} = h^*D_{v+1} . h^*D_{w+1}$
and
$$D_{v+1} . D_{w+1} = \frac{1}{q^{2}} D_v . D_w = \cdots =
\frac{1}{q^{2b}} D_{v+1-b} . D_{w + 1 - b}.$$
On the other hand, $D_i . D_j \in \frac{1}{d} \Z$ with $d$ the determinant
of the intersection matrix for the exceptional divisor of a resolution of $S$.
Thus $D_i . D_{i+1} = D_i^2 = 0$ for $i >>0$. This and the Hodge index theorem
applied to the resolution of $S$, imply that $D_i$ and $D_{i+1}$ are parallel.
So $C_{i}$ and $C_{i+1}$ (and hence $C$ and $f(C)$) are parallel.
Therefore, (1) is true in Case(2). This completes the proof of the lemma.
\end{proof}

\begin{setup}
{\bf Proof of Theorem \ref{Th3-4a} (I)}
\end{setup}

By the assumption, $f^*H_X \sim qH_X$ for an ample line bundle $H_X$.
We will inductively define $\D_i \subset X_i$, $\tau_i : V_i \to X_i$,
$g_{V_i} : V_i \to V_i$, $g_i : X_i \to X_i$,
and big and semi-ample line bundle $H_{V_i}$
with $g_{V_i}^* H_{V_i} \sim qH_{V_i}$.
Define $H_{X_i}$ to be (a large multiple of) the direct image of $H_{V_i}$,
so $g_i^*H_{X_i} \sim qH_{X_i}$ using Lemma \ref{degbir}.
Since $X_i$ is $\Q$-factorial by MMP, $H_{X_i}$ is a big line bundle.
Consider:

\par \vskip 1pc
Property(i): Theorem \ref{Th3-4a} (I) holds for $X_0 \ratmap \cdots
\ratmap X_i$.
$(V_i, \ g_{V_i}, \ X_i \supset \D_i, \ g_i)$ satisfies the four conditions in
Lemma \ref{Mfin}. $H_{V_i}$ is big and semi-ample.
$\dim \D_i \le 2$.

\par \vskip 1pc
The last inequality should follow from the fact:
for a divisorial contraction $\sigma : W \to Z$ between $n$-folds
with exceptional
divisor $E_{W/Z}$, one has $\dim \sigma(E_{W/Z}) \le n-2$;
for a flip $W \ratmap W^+$ with $W \to Z$ and $W^+ \to Z$ the
flipping contractions, one has $\dim E_{W'/Z} \le n-2$ for both
$W' = W, W^+$.

We prove Property(i) ($0 \le i \le r$) by induction.
Set
$$V_0 = X_0, \hskip 1pc \D_0 = \emptyset, \hskip 1pc H_{V_0} := H_X,
\hskip 1pc g_{V_0} = g_0 = f.$$
Then Property(0) holds.
Suppose Property(i) holds for $i \le t$.
If $X_t \to X_{t+1}$ is a divisorial contraction,
then we just apply Lemma \ref{div}.

When $X_t \ratmap X_{t+1} = X_t^+$ is a flip,
we apply Lemma \ref{flip2} and set $\D_{t+1} := \D(X_t^+)$
so that Property(t+1) holds.
Indeed, the first condition in Lemma \ref{flip2}
is satisfied, thanks to Proposition \ref{pol2}. This proves Theorem \ref{Th3-4a} (I).

\begin{setup}
{\bf Proof of Theorem \ref{Th3-4a} (II)}
\end{setup}

By Theorem \ref{descend}, replacing $f$ by its power, we may assume that
$g(C)$ is parallel to $C$ in $N_1(W)$ so that $g : W \to W$ descends to a
finite morphism $h : Y \to Y$; see Lemma \ref{-Kext}.
Set $H_W := H_{X_r}$, a big effective line bundle with $g^*H_W \sim q H_W$.
Now Theorem \ref{Th3-4a} follows from:

\begin{lemma}\label{pi}
$ $
\begin{itemize}
\item[(1)]
$\deg h = q^{\dim Y}$.
\item[(2)]
All eigenvalues of $g_i^* | N^1(X_i)$ and $h^*|N^1(Y)$ are of modulus $q$;
the intersection sheaf $H_Y := I_{V_r/Y}(H_{V_r}^s)$ (with $s = 1 + \dim V_r - \dim Y$)
is a big $\Q$-Cartier integral divisor
such that $h^*H_Y \sim_{\Q} qH_Y$; so $h$ is polarized of degree $q^{\dim Y}$
when $\dim Y \le 2$.
\item[(3)]
If $h$ is polarized, then
$g : W \to W$ is polarized of degree $q^{\dim W}$.
\item[(4)]
Suppose that $h^* | N^1(Y) = q \ \id$.
Replacing
$f$ by its power, we have
$$g_i^* | N^1(X_i) = q \ \id \,\,\,\, (0 \le i \le r).$$
Hence $h$ and $g_i$ are all polarized (see Lemma $\ref{ev}$).
\end{itemize}
\end{lemma}

\begin{proof}
(1) follows from Lemma \ref{ev} and the proof of Lemma \ref{degbir}.

(2) The first part follows from Lemmas \ref{ev} and \ref{degbir}.
We use the birational morphism
$V_r \to X_r = W$ and the big and semi-ample line bundle $H_{V_r}$ in \ref{Th3-4a} (I).
Replacing $H_{V_r}$ by its large multiple, we may assume that $Bs|H_{V_r}| = \emptyset$.
Thus the second part is true
as in Proposition \ref{pol2},
since $I_{V_r/Y}(H_{V_r}^s) = \tau_*(H_{V_r}|V')$,
where $\tau$ is the restriction to $V' := H_1 \cap \cdots \cap H_{s-1}$ of the composite 
$V_r \to W \to Y$,
with $H_i$ general members in $|H_{V_r}|$.
The last part follows from
Theorem \ref{Thev2a} and Lemma \ref{ev}.

(3) We may assume $h^*L \sim qL$ for an ample line bundle
$L$ on $Y$ (using (1)).
The big divisor $H_W$ is $\pi$-ample since
$N_1(W/Y)$ is generated by the class $[C]$.
Thus $H := H_W + t \pi^*L$ is ample for $t >> 0$
(see \cite[Proposition 1.45]{KM}) and $g^*H \sim qH$, so $g$ is polarized.

(4) is true because
$N^1(X_i)$ is spanned by the pullbacks of:
the nef and big divisor $H_{W}$ in \ref{Th3-4a} (I),
the divisors (lying below those divisors in $S(V_j)$, $j \ge i$)
contracted by $X_j \ratmap W$ and
the divisors
in $\pi^*N^1(Y)$, noting that a flip $X_k \ratmap X_{k+1}$ induces
an isomorphism $N^1(X_k) \cong N^1(X_{k+1})$
(see Lemmas \ref{S}, \ref{degbir} and \ref{ev}). This proves Lemma \ref{pi}
and also Theorem \ref{Th3-4a}.
\end{proof}

\begin{setup}
{\bf Proof of Theorem \ref{Thrat3}}
\end{setup}

By Theorem \ref{Th3-4a}, $f$ (replaced by its power) induces a polarized
endomorphism $g : W \to W$ of degree $q^3 > 1$.
Note that $W$ is also rationally connected and
$\Q$-factorial with at worst terminal singularities.
So $K_W$ is not nef.
If the Picard number $\rho(W) = 1$, then
$-K_W$ is ample, and hence $W \cong \P^3$ (so $X$ is rational)
provided that $W$ is smooth, because
every smooth Fano threefold of Picard number one having an endomorphism
of degree $> 1$, is $\P^3$; see \cite{ARV} and \cite{HM}.

Thus, we only need to consider the extremal contraction $\pi: W \to Y$
with $\dim Y = 1, 2$.
Our $Y$ is rational. Note that $\Sing W$ and hence its image
in $Y$ are finite sets, so a general fibre $W_y \subset W$ over $y \in Y$
is smooth.

We apply Theorem \ref{Th3-4a}. Hence each $U \in \{X, W, Y\}$ has an endomorphism $f_U : U \to U$
polarized by an ample line bundle $H_U$ and with $\deg(f_U) = q^{\dim U} > 1$.
Here $f_W = g$ and $f_Y = h$ in notation of Theorem \ref{Th3-4a}.

A polarized endomorphism of degree $> 1$ has a dense set of periodic points
(\cite[Theorem 5.1]{Fa}). Let $y_0$ be a general point with $h(y_0) = y_0$
(after replacing $f$ by its power). Then the fibre $W_0 := W_{y_0} \subset W$
over $y_0 \in Y$
has an endomorphism $g_0 := g|W_0 : W_0 \to W_0$
polarized by the ample line bundle $H_0 := H_W | W_0$ so that
$g_0^*H_0 \sim qH_0$ and $\deg g_0 = q^{\dim W_0} > 1$.
Our $W_0$ is a smooth Fano variety with $\dim W_0 = \dim W - \dim Y$.

Suppose that $\dim Y = 1$. Then $W_0$ is a del Pezzo surface
with a polarized endomorphism of degree $q^2 > 1$. Thus
$K_{W_0}^2 = 6, 8, 9$ (see \cite[Theorem 1.1]{FN} or \cite[Theorem 3]{Zh}; \cite[page 73]{Mi}).
The case $K_{W_0}^2 = 7$ does not occur because $\rho(W/Y) = 1$.
Thus, $W$ (and hence $X$) are rational (see e.g. \cite[\S 2.2]{Is}).

Therefore, we may assume that $\dim Y = 2$. Then $\pi: W \to Y$ is a conic bundle.
$\pi$ is dominated by another conic bundle $\pi': W' \to Y'$
with $W', Y'$ smooth, with $\rho(W'/Y') = 1$ and with
birational morphisms $\sigma_w : W' \to W$ and $\sigma_y: Y' \to Y$
satisfying $\pi \circ \sigma_w = \sigma_y \circ \pi'$ (cf. \cite[the proof of Theorem 4.8]{Mi}).

Let $D'$ be the discriminant of $\pi'$.
If $D' = \emptyset$, then $\pi'$ is a $\P^1$-bundle in the Zariski topology
which is locally trivial for the Brauer group $\Br(Y') = 0$ with
$Y'$ being a smooth projective rational surface, so $W'$ and $X$ are rational.
Thus we may assume that $D' \ne \emptyset$ and $\pi'$ is a standard conic bundle;
see \cite[\S 4.9 and Lemma 4.7]{Mi} for the relevant material.

Let $D$ be the $1$-dimensional part of the discriminant of $\pi$.
Note that $\sigma_{y*}(D') = D$ because every reducible fibre over some $d \in D$
should be underneath only reducible fibres over some $d' \in D'$ and
note that $\sigma_y : Y' \to Y$ is the blowup over the discriminant $D(W/Y)$; see the
construction in \cite[Theorem 4.8]{Mi};
note also that $(\pi')^*E$ is irreducible
for every prime divisor $E \subset X'$
(and especially for those in $D'$).

Our $h : Y \to Y$ satisfies $h^{-1}(D) \subseteq D$
since the reducibility of a fibre $W_d$ over $d$ $\in D$ implies that of
$W_{d'}$ for $d' \in h^{-1}(d)$. So $D \supseteq h^{-1}(D) \supseteq h^{-2}(D)
\supseteq \cdots$.
Considering the number of components, we have $h^{-s}(D) = h^{-s-1}(D)$ for some $s > 1$.
Since $h$ is surjective and applying $h^s$ and $h^{s+1}$, we have
$h^{\pm}(D) = D$. Replacing $f$ by its power, we may assume
$h^{\pm}(D_i) = D_i$ for every irreducible component $D_i$ of $D$.
So $h^*D_i = qD_i$ by Lemma \ref{ev'}.
Hence
$$K_Y + D = h^*(K_Y + D) + G$$
with $G$ an effective Weil divisor. Noting that $h_*H_Y = (\deg(h)/q) H_Y = qH_Y$ and
by the projection formula,
$$H_Y . (K_Y + D) = h_*H_Y . (K_Y + D) + H_Y.G, \,\,\,
(1-q)H_Y.(K_Y + D) = H_Y. G \ge 0.$$
This proves the second assertion below.
For the first, see \cite[Proposition 3.36]{KM} and \cite[Theorem 1.2.7]{MP}.
For the third, see \cite[Lemma 4.1 and Remark 4.2]{Mi}.
The fifth is due to Iskovskikh in his yr 1987 paper in Duke Math. J.
(see e.g. his survey \cite[Theorem 8]{Is}).

\begin{claim}\label{Isk}
$ $
\begin{itemize}
\item[(1)]
$Y$ is $\Q$-factorial with at worst Du Val singularities.
\item[(2)]
If $K_Y+D$ is pseudo-effective, then $K_Y + D \equiv 0$ in $N^1(Y)$.
\item[(3)]
$D'$ is of normal crossing. Every smooth rational component of $D'$
meets at least two points of other components.
\item[(4)]
$\sigma_{y*}(D') = D$.
\item[(5)]
If $\pi'$ is a standard conic bundle, $D'$ is connected and
$D' . F \le 3$ for a free pencil $|F|$ of smooth rational curves,
then $W'$ and hence $W$ and $X$ are rational.
\end{itemize}
\end{claim}

We factor $Y' \to Y$ as $Y' \to \widetilde{Y} \to Y$ with $\widetilde{Y} \to Y$
the minimal resolution. Let $\widetilde{D} \subset \widetilde{Y}$ be the image of $D'$.
Since $D' \ne \emptyset$ and by Claim \ref{Isk} (3) and
the Riemann-Roch theorem, we have $|K_{Y'} + D'| \ne \emptyset$;
the latter implies $K_{\widetilde Y} + \widetilde{D} \sim E$ for some effective divisor.
Hence $K_Y + D \sim \hat{E}$ with $\hat{E} \subset Y$ the image of $E$.
By Claim \ref{Isk} (2), $\hat{E} = 0$ and $K_Y + D \sim 0$.
Thus $\Supp E = \cup_i E_i$ is supported on the exceptional locus of $\widetilde{Y} \to Y$,
so each $E_i$ is a $(-2)$-curve. Now
$h^0(\widetilde{Y}, K_{\widetilde Y} + \widetilde{D}) = 1$.
Our $\widetilde{D}$ is connected and is either a smooth elliptic curve,
or a nodal rational curve,
or a simple loop of smooth rational curves; in fact, one may use Claim \ref{Isk} (3)
and \cite[the proof of Lemma 2.3]{CCZ}.

We assert that $E = 0$. Indeed, since $E$ is negative definite, we may assume that
$E . E_1 < 0$. Then $0 > E_1. (K_{\widetilde Y} + \widetilde{D}) = E_1 . \widetilde{D}$
and hence $E_1 \le \widetilde{D}$.
If $\widetilde{D}$ is irreducible then $E_1 = \widetilde{D}$ and $K_{\widetilde Y}
\sim E - E_1 \ge 0$, contradicting the fact that $\widetilde{Y}$ is a smooth
rational surface. So $\widetilde{D}$ is a simple loop of smooth rational curves
and contains $E_1$.
Thus $0 > E_1 . E_1 + E_1 . (\widetilde{D} - E_1) \ge -2 + 2$ by Claim \ref{Isk} (3).
This is absurd. So our assertion is true and
$K_{\widetilde Y} + \widetilde{D} \sim 0$.

If $\widetilde{Y}$ is ruled with a general fibre $F$
then $\widetilde{D} . F = -K_{\widetilde Y} . F = 2$;
if $\widetilde{Y} = \P^2$, then for a line $F$ we have
$F . \widetilde{D} = 3$.
Denoting by the same $F$ its total transform on $Y'$, we have
$F . D' \le 3$. Thus $W'$ and hence $X$ are rational by Claim \ref{Isk}.
This proves Theorem \ref{Thrat3}.

\begin{setup}
{\bf Proof of Theorem \ref{Thev3}}
\end{setup}

We apply Theorem \ref{Th3-4a}.
By MMP, we may assume that $W$ has no extremal contraction of birational type.
Since $X$ is rationally connected, both $K_X$ and
$K_W$ are non-nef, so there is a contraction $W \to Y$ of an extremal ray.
We have $\dim Y \le 2$. Now Theorem \ref{Thev3} (1) follows from
Theorems \ref{Thev2a} and \ref{Th3-4a} and Lemma \ref{pi} (4) (2).
Indeed, when $\dim Y = 2$, $Y$ is rational with only Du Val singularities
by \cite[Theorem 1.2.7]{MP} and hence $K_Y$ is not trivial in $N^1(Y)$.

Theorem \ref{Thev3} (3) follows from:

\begin{claim}\label{c1Thev3}
Replace $f$ by its power so that $f^*|N^1(X) = q \, \id$. We have:
\begin{itemize}
\item[(1)]
If $M \subset X$ is an irreducible divisor with $\kappa(X, M) = 0$
then $f^*M = q M$.
\item[(2)]
There are only finitely many $f^{-1}$-periodic irreducible divisors $M_i$.
So there is a $v > 0$ such that
$(f^v)^*M_i = q^vM_i$ for all $i$.
The ramification divisor $R_{f^v}$ equals $(q^v-1)\sum_i M_i + \Delta$,
where $\Delta$ is an effective integral divisor containing no any $M_i$.

\item[(3)]
$-K_X \sim_{\Q} \sum_i M_i + \Delta/(q^v -1) \ge 0$ and
$\kappa(X, -K_X) = \kappa(X, \sum M_i - K_X) \ge 0$.
\end{itemize}
\end{claim}

\begin{proof}
Since $q(X) = 0$, we have $f^*M \sim_{\Q} qM$ for every
irreducible integral divisor $M$.
Suppose that $\kappa(X, M) = 0$. Since $f^*M \sim_{\Q} qM$, we have $f^{-1}(M) = M$.
Then (1) follows.

Suppose that $M_i$ ($1 \le i \le N$) are $f^{-1}$-periodic,
so a power $h_N = f^{s(N)}$ of $f$ satisfies $h_N^{-1}(M_i) = M_i$ for all $1 \le i \le N$.
Then $h_N^*M_i = q^{s(N)}M_i$ and
$K_X + \sum M_i = h_N^*(K_X + \sum M_i) + \Delta_N \sim_{\Q} q^{s(N)}(K_X + \sum M_i) + \Delta_N$,
where $\Delta_N$ is an effective
integral divisor containing no any $M_i$.
Thus $-K_X \sim_{\Q} \sum_{i=1}^N M_i + \Delta/(q^{s(N)} -1) \ge 0$, which also implies (3).
Multiplying the above equivalence by $\dim X - 1 = 2$ copies of an ample divisor $H$,
we see that $N$ is bounded. This proves (2).
\end{proof}

We now prove Theorem \ref{Thev3} (2).
By Theorem \ref{Thrat3}, we may assume that the end product of MMP for $X$
is of Picard number one, i.e., there is a composite
$X = X_0 \ratmap X_1 \cdots \ratmap X_r$ of divisorial contractions and flips
such that $\rho(X_r) = 1$, so $-K_{X_r}$ is ample because all $X_i$ are rationally
connected with only $\Q$-factorial terminal singularities by MMP.
Let $g_i : X_i \ratmap X_i$ be the dominant rational map induced from $f : X \to X$
(with $g_0 = f$).

\begin{claim}\label{c2Thev3}
Replacing $f$ by its positive power, we have:
\begin{itemize}
\item[(1)]
For all $0 \le t \le r$, our $g_t$ is holomorphic with $g_t^*|N^1(X_t) = q \, \id$.
Let $E_t' \subset X_t$ be zero (resp. the (irreducible) exceptional divisor)
when $X_t \ratmap X_{t+1}$ is a flip (resp. $X_t \to X_{t+1}$ is divisorial).
Then the strict transform $E_t \subset X$ of $E_t'$
satisfies $f^{-1}(E_t) = E_t$.
\item[(2)]
$N^1(X)$ is spanned by $K_X$ and those $E_t$ in (1).
Let $E = \sum E_t$.
\end{itemize}
\end{claim}

\begin{proof}
(1) can be proved by ascending induction on the index $t$ of $X_t$.
Suppose (1) is true for $t$.
Since $g_t^*$ is scalar, we may assume that both $g_t^{\pm}$ preserve
the extremal ray corresponding to the birational map $X_t \ratmap X_{t+1}$,
so $g_t$ descends to the holomorphic $g_{t+1}$ as in the proof of Theorem
\ref{Th3-4a}, and also the last part of (1) is true.
The scalarity of $g_t^*$ implies that of $g_{t+1}^*$ because
$N^1(X_{t+1})$ is isomorphic to (resp. regarded as a subspace of) $N^1(X_t)$ via the pullback
when $X_t \ratmap X_{t+1}$ is a flip (resp. $X_t \to X_{t+1}$ is divisorial);
see \cite[the proof of Proposition 3.37]{KM}.

(2) is true because $N^1(X_r)$ is generated by $K_{X_r}$,
$N^1(X_t)$ is isomorphic to $N^1(X_{t+1})$ (resp. spanned by $E_t'$
and the pullback of $N^1(X_{t+1}$)) when $X_t \ratmap X_{t+1}$ is
a flip (resp. divisorial).
\end{proof}

To conclude Theorem \ref{Thev3} (2),
take an ample divisor $H \subset X$. By Claim \ref{c2Thev3},
we can write $H \sim_{\Q} \sum a_t E_t + b(-K_X)$.
So $H \le m(E - K_X)$ for some $m \ge 1$, since $\kappa(X, -K_X) \ge 0$.
This and Claim \ref{c1Thev3} (3) and Claim \ref{c2Thev3} (1)
imply $\kappa(X, -K_X) = \kappa(X, E - K_X) \ge \kappa(X, H) = \dim X$.
Thus, $-K_X$ is big. Theorem \ref{Thev3} (2) is proved.

\begin{setup}
{\bf Proof of Theorem \ref{Corrat3}}
\end{setup}

Since $X$ is Fano,
$X$ is rationally connected (by Campana and Koll\'ar-Miyaoka-Mori),
and $\overline{\NE}(X)$ has only finitely
many extremal rays all of which are $K_X$-negative (cf. \cite[Theorem 3.7]{KM}).
Let $X \to X_1$ be the smooth blowdown such that $X_1$ is a
primitive (smooth) Fano threefold in the sense of \cite{MM}.
If $\rho(X) \ge 2$,
by \cite[Theorem 5]{MM}, $X_1$ has an extremal contraction of conic bundle type.
Now Theorem \ref{Corrat3} follows from Theorem \ref{Thrat3}.
%
%
%
%
%
%

\begin{setup}
{\bf Proof of Theorem \ref{FinExt1}}
\end{setup}
By Lemma \ref{S}, we may assume that $S(X)$ is a finite set.
We may also assume $\rho(X) \ge 3$.
Suppose that $R_i := \R_{\ge 0}[C_i]$ ($i \ge 1$) are pairwise distinct
$K_X$-negative extremal rays with $\pi_i : X \to Y_i$ the corresponding
contraction each of which is either divisorial or of Fano type (i.e., $\dim Y_i \le 2$).
We can take the generator $C_i$ to be an irreducible curve in the fibre of $\pi_i$.
Since $3 \le \rho(X) = \rho(Y_i) + 1$, we have $\rho(Y_i) \ge 2$ and hence $\dim Y_i \in \{2, 3\}$.

If $\pi_i$ is divisorial, we let $E_i$ be the exceptional divisor of $\pi_i$;
then $E_i$ is necessarily irreducible and is in the finite set $S(X)$.
If $\pi_i$ is of Fano type (and hence onto a surface $Y_i$), then
$Y_i$ is a rational surface with at worst Du Val singularities (cf. \cite[Theorem 1.2.7]{MP});
for each $G \in S(Y_i)$, the divisor $\pi_i^*G$ is irreducible and in $S(X)$.

The claim below follows from the fact that $\rho(X/Y_i) = 1$.

\begin{claim}\label{res}
Suppose that either $D$ is the exceptional divisor $E_i$ for a divisorial contraction $\pi_i : X \to Y_i$,
or $D = \pi_i^*G$ for a Fano contraction $\pi_i : X \to Y_i$ to a surface with $G \subset Y_i$ an irreducible curve.
Then $N^1(X) | D$, as a subspace of $N_1(X)$, is of rank $\le 2$ and contains the extremal ray
$R_i$ of $\overline{\NE}(X)$.
\end{claim}

Suppose, after replacing with an infinite subsequence, that each $\pi_i$ is either divisorial
and we let $D_i := E_i$, or is of Fano type with $S(Y_i) \ne \emptyset$ and we let
$D_i = \pi_i^*G$ for some $G \in S(Y_i)$. Since $D_i \in S(X)$ and $S(X)$ is finite,
we may assume that $D_1 = D_2 = \cdots$ after replacing with an infinite subsequence.
If $N^1(X) | D_i \subset N_1(X)$ contains only one extremal ray, i.e., $R_i$, then $R_1 = R_2$, absurd.
If $N^1(X) | D_i$ has two extremal rays $R_i, R_i'$, then either $R_i = R_j$ for some $i \ne j$ absurd;
or $R_2 = R_1' = R_3$, absurd again.

Thus, replacing with an infinite subsequence, we may assume that for every $i
\ge 1$, $\pi_i$ is of Fano type and $S(Y_i) = \emptyset$. 
Hence $Y_i$ is relatively minimal, $\rho(Y_i) = 2$
and there is a $\P^1$-fibration $Y_i \to B_i \cong \P^1$ with every fibre irreducible, noting that
$K_{Y_i}$ is not pseudo-effective (cf. \cite[Theorem 3.2]{Sa}).
Take a general fibre $X_{b_i}$ of the composite $X \to Y_i \to B_i$
which is a smooth relatively minimal ruled surface, 
noting that $\Sing X$ and hence its image in $B_i$ are finite sets.
Then $R_i . X_{b_i} = 0$.

Now $\rho(X) = \rho(Y_i) + 1 = 3$. Any three of $C_i$ are linearly independent in $N_1(X)$
and hence form a basis;
otherwise, $C_3 = a_1C_1 + a_2C_2$ say with $a_1 > 0, a_2 \ge 0$ and hence $R_1 = R_3$,
since $R_3$ is extremal. This is a contradiction.

Suppose that $R_1 . X_{b_i} = 0$, i.e., $\pi_1(X_{b_i}) \ne Y_1$ for $i = 2, 3, 4$.
Then $X_{b_i} = \pi_1^*M_i$ for an irreducible curve $M_i \subset Y_1$ since $\rho(X/Y_1) = 1$.
Since $\rho(Y_1) = 2$ and $q(Y_1) = 0$, we may assume that $M_4 \sim_{\Q} a_2 M_2 + a_3M_3$
and hence $X_{b_4} \sim_{\Q} a_2X_{b_2} + a_3X_{b_3}$. Note that $0 = X_{b_4}^2 = 2a_2a_3 X_{b_2} X_{b_3}$.
After relabeling, we may assume that $X_{b_3}$ and $X_{b_4}$ are parallel in $N^1(X)$.
Then $X_{b_3} = \pi_1^*M_3$ is perpendicular to all of $C_1, C_3, C_4$, a basis of $N_1(X)$.
So $X_{b_3} = 0$ in $N^1(X)$. This is absurd.

Therefore, we may assume that $\pi_1(X_{b_i}) = Y_1$ for all $i \ge 2$, after replacing with a subsequence.
Since $S(Y_1) = \emptyset$ and $\rho(Y_1) = 2$, our
$\overline{\NE}(Y_1)$ is generated by two extremal pseudo-effective divisors $L_1, L_1'$
with $L_1^2 = (L_1')^2 = 0$. We may assume that $L_1$ is a fibre of $Y_1 \to B_1$.
Let $M_i, M_i' \in N_1(X_{b_i})$ (which are necessarily linearly independent
and hence form its basis)
be respectively the pullbacks of $L_1, L_1'$, via $\pi_1|X_{b_i}$.
We may assume that $C_i$ (a fibre of $\pi_i$) belongs to $N_1(X_{b_i})$.
Then $C_i = e M_i + e'M_i'$ in $N_1(X_{b_i})$. Since
$0 = C_i^2 = 2ee' M_i . M_i'$  on $X_{b_i}$  and $M_i . M_i' = \deg(\pi_1|X_{b_i}) L_1 . L_1' > 0$,
we have $ee' = 0$. So $C_i$ is parallel to $M_i$ or $M_i'$ in $N_1(X)$.
If $C_i$ is parallel to $M_i = X_{b_i} \cap \pi_1^*L_1$ for $i = r, s, t$
then by Claim \ref{res} applied to $N^1(X)|\pi_1^*L_1$, two of the (extremal) $C_i$
are parallel to each other in $N_1(X)$, contradicting the fact that $R_i$'s are all distinct.
If $C_i$ is parallel to $M_i'$ for $i = u, v, w$, then $(\pi_1|X_{b_i})_*C_i$ is parallel to
$L_1'$ and we may assume that $L_1'$ is an irreducible curve. Applying Claim \ref{res} to $N^1(X) | \pi_1^*L_1'$,
we get a similar contradiction.
This proves Theorem \ref{FinExt1}.

\end{large}
%
%

%
%
%
\par \vskip 2pc \noindent
De-Qi Zhang
\par \noindent
Department of Mathematics, National University of Singapore
\par \noindent
2 Science Drive 2, Singapore 117543, SINGAPORE
\par \noindent
E-mail: matzdq@nus.edu.sg

\newpage

\par \vskip 2pc

$$
\text{APPENDIX}
$$
%
%
%
%




\setlength{\topmargin}{0cm}
\setlength{\oddsidemargin}{0cm}
\setlength{\evensidemargin}{0cm}
\setlength{\marginparwidth}{0cm}
\setlength{\marginparsep}{0cm}
\setlength{\headheight}{12.5pt}
\setlength{\headsep}{25pt}
\setlength{\footskip}{30pt}
\renewcommand{\baselinestretch}{1.2}

\pagestyle{headings}
\theoremstyle{plain}
    \newtheorem{thm}{Theorem}
    \newtheorem{lem}[thm]{Lemma}%
\theoremstyle{remark}
    \newtheorem{rem}[thm]{Remark}
    \newtheorem{fact}[thm]{Fact}
    \newtheorem{nota}[thm]{Notation}

\renewcommand{\theequation}{*\arabic{equation}}

\newcommand{\BZZ}{\mathbb{Z}}
\newcommand{\BPP}{\mathbb{P}}
\newcommand{\BQQ}{\mathbb{Q}}
\newcommand{\BRR}{\mathbb{R}}
\newcommand{\BCC}{\mathbb{C}}
\newcommand{\SK}{\mathcal{K}}
\newcommand{\ep}{\varepsilon}
\newcommand{\class}{\operatorname{cl}}


\title[]{Termination of extremal rays of fibration type
for the iteration of surjective endomorphisms}
\author[]{Yoshio Fujimoto}
\address[Yoshio Fujimoto]{
Department of Mathematics,
Faculty of Medical Science, \endgraf
Nara Medical University,
840 Shijo-cho Kashihara-city Nara 634-8521 Japan.}
\email{yoshiofu@naramed-u.ac.jp}

\author[]{Noboru Nakayama}
\address[Noboru Nakayama]{%
Research Institute for Mathematical Sciences\endgraf
Kyoto University, Kyoto 606-8502 Japan}
\email{nakayama@kurims.kyoto-u.ac.jp}

\maketitle

The purpose of this note is to prove the following:

\begin{thm}\label{thm1}
Let \( X \) be a normal projective variety defined over an algebraically
closed field of characteristic zero such that \( X \)
has only log-terminal singularities.
Let \( R \subset \NEb(X) \) be an extremal ray such that \( K_{X}R < 0 \)
and the associated contraction morphism \( \Cont_{R} \) is a fibration to
a lower-dimensional variety.
Then, for any surjective endomorphism \( f \colon X \to X \),
there exists a positive integer \( k \) such that
\( (f^{k})_{*}(R) = R \) for the automorphism
\( (f^{k})_{*} \colon \NN_{1}(X) \xrightarrow{\isom} \NN_{1}(X)\)
induced from the iteration \( f^{k} = f \circ \cdots \circ f \).
\end{thm}

A special case is proved in Theorem~2.13 of
a recent paper \cite{Z} of D.-Q.~Zhang.
We extend and simplify the idea of Zhang.
The authors express their gratitude to Professor De-Qi Zhang
for informing his paper \cite{Z}.

\begin{nota}
For a normal projective variety \( X \), let \( \NN^{1}(X) \) denote the vector space
\( \operatorname{NS}(X) \otimes \BRR\) for the N\'eron-Severi group
\( \operatorname{NS}(X) \). The dimension of \( \NN^{1}(X) \)
is called the \emph{Picard number}
and is denoted by \( \rho(X) \).
The numerical equivalence class \( \class(D) \) of
a Cartier divisor \( D \) on \( X \) is regarded as an element of
\( \NN^{1}(X) \).
The dual vector space of \( \NN^{1}(X) \) is denoted by \( \NN_{1}(X) \), i.e.,
\( \NN_{1}(X) = \Hom(\operatorname{NS}(X), \BRR) \).
An element \( u \in \NN^{1}(X) \) is regarded as a linear function
on \( \NN_{1}(X) \). We denote by \( u^{\perp} \)
the kernel of \( u \colon \NN_{1}(X) \to \BRR \).
The cone \( \NE(X) \) of the numerical equivalence classes \( \class(Z) \)
of the effective \( 1 \)-cycles \( Z \) on \( X \)
is defined in \( \NN_{1}(X) \),
by the intersection pairing \( D \mapsto DZ \in \BZZ \)
for Cartier divisors \( D \) on \( X \).

The closure of \( \NE(X) \) in \( \NN_{1}(X) \)
is denoted by \( \NEb(X) \), which is a \emph{strictly convex} cone, i.e.,
\( \NEb(X) + \NEb(X) \subset \NEb(X) \) and \( \NEb(X) \cap (-\NEb(X)) = \{0\} \).
An \emph{extremal ray} \( R \) of \( \NEb(X) \) is
by definition a one-dimensional face of the cone \( \NEb(X) \), i.e.,
\( R = \BRR_{\geq 0} v = u^{\perp} \cap \NEb(X)\)
for some \( 0 \ne v \in \NEb(X) \) and for some \( u \in \NN^{1}(X) \)
which is non-negative on \( \NEb(X) \) as a function on \( \NN_{1}(X) \).
For a Cartier divisor \( D \) on \( X \), \( DR > 0 \) means that
the functional \( \class(D) \) on \( \NN_{1}(X) \) is
positive on \( R \setminus \{0\} \).
The meanings of \( DR = 0 \) and \( DR < 0 \) are similar.
\end{nota}

\begin{fact}[\cite{Ka}]\label{fact:cont}
Let \( X \) be a normal projective variety with only log-terminal singularities,
i.e., \( (X, 0) \) has only log-terminal singularities in the sense of \cite{Ka}.
For an extremal ray \( R \) of \( \NEb(X) \) with \( K_{X}R < 0 \),
there exist a proper surjective morphism \( \Cont_{R} \colon X \to Y \)
onto a normal projective variety \( Y \) satisfying the following two conditions:
\begin{enumerate}
\item \label{fact:cont:1} Every fiber of \( \Cont_{R} \) is connected.

\item \label{fact:cont:2} For an irreducible closed curve \( C \) on \( X \),
\( \Cont_{R}(C) \) is a point if and only if  \( \class(C) \in R\).
\end{enumerate}
The morphism \( \Cont_{R} \) is uniquely determined by the conditions
\eqref{fact:cont:1} and \eqref{fact:cont:2}, and is called
the \emph{contraction morphism} associated with \( R \).
The following property holds by \cite[Corollary~4.4]{Ka}:
\begin{enumerate}
    \addtocounter{enumi}{2}
\item \label{fact:cont:3}
If \( D \) is a Cartier divisor on \( X \) with \( DR = 0 \), then
\( D \sim \Cont_{R}^{*}(E) \) for a Cartier divisor \( E \) on \( Y \).
\end{enumerate}
\end{fact}

\begin{rem}\label{rem:rho=rho}
Let \( f \colon X \to Y \) be a surjective morphism between
normal projective varieties.
Then, we have
the pullback homomorphism \( f^{*} \colon \NN^{1}(Y) \to \NN^{1}(X)\) which is
well-defined by \( f^{*}(\class(D)) := \class(f^{*}(D))  \) for Cartier divisors
\( D \) on \( Y \).
We have also the push-forward homomorphism
\( f_{*} \colon \NN_{1}(X) \to \NN_{1}(Y) \) as the dual of \( f^{*} \).
Here, for any irreducible closed curve \( C \) on \( X \), we have
\( f_{*}(\class(C)) = \class(f_{*}(C))\) for the \( 1 \)-cycle
\[ f_{*}(C) = \begin{cases}
\deg (C/f(C)) C, & \text{ if \( f(C) \) is not a point}; \\
0, & \text{ otherwise.}
\end{cases}\]
Since \( f \) is surjective,
\( f^{*} \colon \NN^{1}(Y) \to \NN^{1}(X) \) is injective and
\( f_{*} \colon \NN_{1}(X) \to \NN_{1}(Y) \) is surjective.
Assume that \( \rho(X) = \rho(Y) \). Then
\( f^{*} \) and \( f_{*} \) above are both isomorphisms, since
\( \NN^{1}(X) \) and \( \NN^{1}(Y) \) have the same dimension.
In particular, we have \( f_{*}(\NEb(X)) = \NEb(Y) \) from the
obvious equality \( f_{*}(\NE(X)) = \NE(Y) \).
Moreover, \( f \) is a finite morphism;
in fact, \( f(C) \) is not a point for
any irreducible closed curve \( C \) on \( X \) by \( f_{*}(\class(C)) \ne 0 \).
\end{rem}

\begin{lem}\label{lem:ray}
In the situation of Theorem~\ref{thm1},
\( f_{*}(R) \) is also an extremal ray of \( \NEb(X) \) such that
\( K_{X}f_{*}(R) < 0 \).
\end{lem}

\begin{proof}
The push-forward map \( f_{*} \colon \NN_{1}(X) \to \NN_{1}(X) \) is an automorphism
preserving the cone \( \NEb(X) \). Thus, \( f_{*}(R) \) is extremal.
Let \( E_{f} \) be the ramification divisor of \( f \colon X \to X \), i.e.,
\( K_{X} = f^{*}(K_{X}) + E_{f} \).
Since \( E_{f} \) is effective, the restriction of \( E_{f} \) to a general fiber of
\( \Cont_{R} \) is also effective. Hence,
\( E_{f}\gamma \geq 0 \) for a general curve
\( \gamma \) contracted to a point by \( \Cont_{R} \).
Thus \( 0 > K_{X}\gamma \geq (f^{*}K_{X}) \gamma = K_{X} (f_{*}\gamma) \).
Therefore, \( K_{X}f_{*}(R) < 0 \).
\end{proof}

\begin{nota}
For the extremal ray \( R \) in Theorem~\ref{thm1},
let \( R_{k} \) be the extremal ray \( f^{k}_{*}(R) \) for \( k \geq 0 \).
By Fact~\ref{fact:cont} and Lemma~\ref{lem:ray}, we have the associated
contraction morphism \( \Cont_{R_{k}} \), which is denoted by
\( \pi_{k} \colon X \to Y_{k} \).
Then, \( \pi_{k+1} \circ f = h_{k} \circ \pi_{k}\)
for a finite surjective morphism \( h_{k} \colon Y_{k} \to Y_{k+1} \)
by the condition \eqref{fact:cont:2} in Fact~\ref{fact:cont};
in particular, we have the following commutative diagram:
\[ \begin{CD}
X @= X @>{f}>> X @>{f}>>
\cdots @>{f}>> X @>{f}>> X @>{f}>> \cdots \\
@V{\pi}VV @V{\pi_{0}}VV @V{\pi_{1}}VV @. @V{\pi_{k}}VV @V{\pi_{k+1}}VV  \\
Y @= Y_{0} @>{h_{0}}>> Y_{1} @>{h_{1}}>>
\cdots @>>> Y_{k} @>{h_{k}}>> Y_{k+1} @>{h_{k+1}}>> \cdots
\end{CD}\]
Here, we simply write \( \pi = \pi_{0} \) and \( Y = Y_{0} \).
We define \( m := \dim Y\) and
\( \rho := \rho(X) - 1 \geq 0\). Then \( m = \dim Y_{k} \),
\( \rho = \rho(Y_{k}) \), and
\( h_{k}^{*} \colon \NN^{1}(Y_{k+1}) \to \NN^{1}(Y_{k}) \)
is an isomorphism for any \( k \geq 0\).
\end{nota}

\begin{lem}\label{lem:rho}
\emph{Theorem~\ref{thm1}} is true if \( \rho \leq 1 \).
\end{lem}

\begin{proof}
Assume that \( \rho = \rho(X) - 1 = 0\). Then \( \NN_{1}(X) \) is one-dimensional and
\( \NEb(X) \) is just a single ray. Thus \( R_{k} = R \) for any \( k \).
Assume next that \( \rho = \rho(X) - 1 = 1 \).
Then \( \NEb(X) \) has exactly two extremal rays.
Hence, \( f^{2}_{*} \) preserves each extremal ray.
Therefore, \( R = R_{2k} \) for any \( k \).
\end{proof}

\begin{lem}\label{lem:Dm=0}
Let \( D \) be a Cartier divisor on \( Y \) such that
\( \pi^{*}(D) R_{k} = 0 \) for some \( k \geq 1 \).
If the self-intersection number \( D^{m} \ne 0 \), then \( R = R_{k} \).
\end{lem}

\begin{proof}
By the property \eqref{fact:cont:3} in Fact~\ref{fact:cont}
of the contraction morphism of an extremal ray, we have
a Cartier divisor \( D_{k} \) on \( Y_{k} \) such that
\( \pi^{*}(D) \sim \pi_{k}^{*}(D_{k}) \).
Let \( A \) be an ample divisor on \( X \). Then the product
\( \pi^{*}(D)^{m}A^{n - m - 1} \) in the Chow ring of \( X \) is
numerically equivalent to
\( \delta Z \) for a non-zero effective \( 1 \)-cycle \( Z \) and for
\( \delta := D^{m} \ne 0 \).
Thus,
\[ \pi^{*}(L)Z = \delta^{-1}\pi^{*}(LD^{m})A^{n - m - 1} = 0
\quad \text{and} \quad
\pi_{k}^{*}(L_{k}) Z = \delta^{-1} \pi_{k}^{*}(L_{k} D_{k}^{m})A^{n - m - 1} = 0 \]
for any Cartier divisor \( L \) on \( Y \) and
any Cartier divisor \( L_{k} \) on \( Y_{k} \).
In particular,
the numerical equivalence class \( \class(Z) \) is contained in \( R \cap R_{k}\).
Therefore, \( R = R_{k} \).
\end{proof}

\begin{proof}[Proof of \emph{Theorem~\ref{thm1}}]
We shall derive a contradiction from the converse assumption
that \( R \ne R_{k} \) for any \( k \geq 1 \).
Then, \( R_k \ne R_j\) for any \(j \ne k\),
since \(f_* \colon \NN_1(X) \to \NN_1(X)\)
is an automorphism by Remark~\ref{rem:rho=rho}.
We have \( \rho \geq 2\) by Lemma~\ref{lem:rho}.
In particular, \( \dim Y = m \geq 2 \).
Let \( \{H_{1}, \ldots, H_{\rho}\} \) be a set of ample divisors of \( Y \)
such that \( \{\class(H_{1}), \ldots, \class(H_{\rho})\} \)
is a basis of \( \NN^{1}(X) \).
We have \( (\pi^{*}H_{i}) R_{k} > 0 \) for any \( 1 \leq i \leq \rho\) and
\( k \geq 1 \) by the property \eqref{fact:cont:3} in Fact~\ref{fact:cont},
since \( R \ne R_{k} \).
Hence, we can define a positive rational number \( a^{(j)}_{k} \)
for \( 2 \leq j \leq \rho\) and \( k \geq 1 \) by the equation:
\begin{equation}\label{eq:ajk}
\pi^{*}(H_{j} - a^{(j)}_{k}H_{1}) \cdot R_{k} = 0.
\end{equation}
Then \( (H_{j} - a^{(j)}_{k}H_{1})^{m} = 0 \) for any \( j \) and \( k \)
by Lemma~\ref{lem:Dm=0}.
On the other hand, for each \( 2 \leq j \leq \rho\),
there exist at most \( m \) solutions for \( x \in \BCC \) of the equation:
\( (H_{j} - xH_{1})^{m} = 0 \).
Then, there exist rational numbers \( \alpha_{2} \), \ldots, \( \alpha_{\rho} \)
such that, for infinitely many integers \( k \), the equalities
\( \alpha_{j} = a^{(j)}_{k} \) hold for any \( 2 \leq j \leq \rho \).
In fact, we can find a rational number \( \alpha_{2} \) such that the set \( S_{2} \)
of positive integers \( k \) with
\( \alpha_{2} = a^{(2)}_{k} \) is infinite.
Next, we can find a rational number \( \alpha_{3} \) such that the set \( S_{3} \)
of integers \( k \in S_{2} \) with \( \alpha_{3} = a^{(3)}_{k} \) is infinite.
If the rational numbers
\( \alpha_{j} \) with the sets \( S_{j} \) up to \( l < \rho\) are selected,
then we can find a rational number \( \alpha_{l+1} \) such that the set
\( S_{l+1} \) of integers \( k \in S_{l} \) with \( \alpha_{l+1} = a^{(l+1)}_{k} \)
is infinite.
In this way, we can find \( \alpha_{2} \), \( \alpha_{3} \),
\ldots, \( \alpha_{\rho} \) satisfying the required property.

The real vector subspace
\[ F := \pi^{*}(\class(H_{2} - \alpha_{2}H_{1}))^{\perp} \cap \cdots \cap
\pi^{*}(\class(H_{\rho} - \alpha_{\rho}H_{1}))^{\perp} \subset \NN_{1}(X)\]
is two-dimensional, since \( \pi^{*}(\class(H_{2} - \alpha_{2}H_{1})) \), \ldots,
\( \pi^{*}(\class(H_{\rho} - \alpha_{\rho}H_{1})) \) are linearly independent.
We have \( R_{k} \subset F  \) for infinitely many \( k \) by the choice of
\( \alpha_{2} \), \ldots, \( \alpha_{\rho} \) and by \eqref{eq:ajk}.
This is a contradiction, since there exist at most two extremal rays
of \( \NEb(X) \) contained in the two-dimensional vector subspace \( F \).
Thus, we are done.
\end{proof}

\begin{rem}\label{rem:Th1}
In Theorem~\ref{thm1}, we can not allow the case where
\( \Cont_{R} \) is a birational morphism.
In fact, there exist a smooth projective surface \( X \) with an automorphism \( f \)
and a \( (-1) \)-curve \( \gamma \) on \( X \)
such that \( \{f^{k}(\gamma) \mid k \geq 0\} \) is infinite.
Here, \( R = \BRR_{\geq 0}\class(\gamma) \) is an extremal ray with \( K_{X}R < 0 \)
and  \( f^{k}_{*}(R) = \BRR_{\geq 0}\class(f^{k}(\gamma)) \) for the \( (-1) \)-curve
\( f^{k}(\gamma) \). Thus \( f^{k}_{*}(R) \ne R \) for any \( k \).
One of such a surface \( X \) is given as a
blown up surface of \( \BPP^{2} \)
whose center is the intersection of two sufficiently general cubic curves.
In fact, \( X \) is a rational elliptic surface and any exceptional curve of the blowing up
is a section of the elliptic fibration.
Let \( \Gamma_{0} \) and \( \Gamma_{1} \) be two exceptional curves.
Let \( X_{K} \) be the generic fiber of the elliptic fibration and
\( P_{i} \) the point \( \Gamma_{i}|_{X_{K}} \)
defined over the function field \( K \) of the base curve.
We give a group structure of the elliptic curve \( X_{K} \) such that \( P_{0} \)
is the zero element. Then, \( P_{1} \) is not torsion by
the choice of cubic curves.
The translation mapping \( X_{K} \to X_{K} \) by \( P_{1} \) gives rise to
a birational automorphism \( f \colon X \to X \), which is in fact regular,
since the elliptic surface \( X \) is relatively minimal over the base curve.
Therefore, \( f \) is an automorphism of infinite order and
\( f^{k}(\Gamma_{1}) \ne \Gamma_{1} \) for any \( k \).
Thus, the conclusion of Theorem~\ref{thm1} does not hold for
\( X \), \( f \), and \( R = \BRR_{\geq}\class(\Gamma_{1}) \).
\end{rem}

\end{document}